\documentclass[11pt]{article}
\usepackage{graphicx,fullpage,latexsym,amssymb,amsmath,url}
\usepackage{epstopdf}
\usepackage{enumitem,color,cite}
\usepackage{booktabs}
\usepackage{mathdots}

\newcommand{\BEAS}{\begin{eqnarray*}}
\newcommand{\EEAS}{\end{eqnarray*}}
\newcommand{\BEA}{\begin{eqnarray}}
\newcommand{\EEA}{\end{eqnarray}}
\newcommand{\BEQ}{\begin{equation}}
\newcommand{\EEQ}{\end{equation}}
\newcommand{\BIT}{\begin{itemize}}
\newcommand{\EIT}{\end{itemize}}
\newcommand{\eg}{{\it e.g.}}
\newcommand{\ie}{{\it i.e.}}
\newcommand{\ones}{\mathbf 1}
\newcommand{\reals}{{\mbox{\bf R}}}
\newcommand{\complex}{{\mbox{\bf C}}}
\newcommand{\Rank}{\mathop{\bf rank}}
\newcommand{\Tr}{\mathop{\bf tr}}
\newcommand{\diag}{\mathop{\bf diag}}
\newcommand{\Expect}{\mathop{\bf E{}}}

\newtheorem{theorem}{Theorem}
\newtheorem{lemma}{Lemma}
\newcommand{\qed}{\hfill$\Box$}
\newcommand{\proof}{\noindent\emph{Proof.}\ }

\newcommand{\herm}{\mathbf{H}}
\renewcommand{\j}{\mathrm{j}}
\DeclareMathOperator{\newRe}{\mathrm{Re}}
\renewcommand{\Re}{\newRe}
\DeclareMathOperator{\newIm}{\mathrm{Im}}
\renewcommand{\Im}{\newIm}

\DeclareMathOperator{\conv}{\mathrm{conv}}

\title{Semidefinite representations of gauge functions
for structured low-rank matrix decomposition}
\author{Hsiao-Han Chao\thanks{
    Electrical Engineering Department, University of California, 
    Los Angeles. Email:
    {\tt suzihchao@ucla.edu}, {\tt vandenbe@ucla.edu}.
    Research partially supported by NSF Grants 
    1128817 and 1509789.}
   \and 
   Lieven Vandenberghe\footnotemark[1]}
\date{}

\begin{document}
\maketitle
\begin{abstract}
This paper presents generalizations of semidefinite programming 
formulations of 1-norm optimization problems over infinite 
dictionaries of vectors of complex exponentials, which were 
recently proposed for superresolution, gridless compressed sensing, 
and other applications in signal processing.
Results related to the generalized Kalman-Yakubovich-Popov 
lemma in linear system theory provide simple, constructive proofs of the 
semidefinite representations of the penalty functions used in these 
applications.  The connection leads to several extensions to 
gauge functions and atomic norms for sets of vectors 
parameterized via the nullspace of matrix pencils.  
The techniques are illustrated with examples of
low-rank matrix approximation problems arising in spectral estimation 
and array processing.
\end{abstract}

\section{Introduction} \label{s-intro}
The notion of atomic norm introduced in \cite{CRPW:12} gives a 
unified description of convex penalty functions that extend 
the $\ell_1$-norm penalty, used to promote sparsity in 
the solution of an optimization problem,
to various other types of structure. 
The atomic norm associated with a non-empty set $C$ is defined as
the gauge of its convex hull, \ie, the convex function
\BEA
g(x) 
& = & \inf {\{t \geq 0 \mid x \in t \conv{C}\}} \nonumber \\
&= & \inf {\{\sum_{k=1}^r \theta_k \mid 
x =\sum_{k=1}^r \theta_k a_k, \; \theta_k \geq 0, \; a_k\in C\}}.
\label{e-atomic}
\EEA
This function is convex, nonnegative, positively homogeneous, 
and zero if $x=0$.  It is not necessarily a norm,
but it is common to use the term `atomic norm' even when $g$ 
is not a norm.
When used as a regularization term in an optimization problem, 
the function $g(x)$ defined in~(\ref{e-atomic}) promotes the property 
that $x$ can be expressed as 
a nonnegative linear combination of a small number of elements 
(or `atoms') of $C$.

The best known examples of atomic norms are the vector $\ell_1$-norm 
and the matrix trace norm.
The $\ell_1$-norm of a real or complex $n$-vector is the atomic norm
associated with $C = \{ se_k \mid |s|=1, \; k=1,\ldots,n\}$, 
where $e_k$ is the $k$th unit vector of length $n$. 
The matrix trace norm (or nuclear norm) 
is the atomic norm for the set of rank-1 matrices with unit 
norm.  Specifically, the trace norm on $\complex^{n\times m}$ 
is the atomic norm for $C=\{vw^H \mid \|v\| =\|w\|=1\}$,   
where $w^H$ is the conjugate transpose and $\|\cdot\|$ denotes
the Euclidean norm.
Many other examples are discussed in \cite{CRPW:12,BTR12,TBSR:13}. 

The atomic norm associated with the set 
\BEQ \label{e-Ce}
 C_\mathrm{e} 
 = \{ \gamma\, (1, e^{\j\omega}, \ldots, e^{\j (n-1) \omega}) 
 \in \complex^n \mid \omega \in [0,2\pi), \; |\gamma|=1/\sqrt n\},
\EEQ
where $\j = \sqrt{-1}$, has been studied extensively in recent 
research in signal processing
\cite{dCG12,dCGH15,CaF:14,BTR12,TBSR:13,YX14,LC14,YX15,CC15}.
It is known that the atomic norm for this set is the optimal value of 
the semidefinite program (SDP)
\BEQ \label{e-Ce-sdp}
 \begin{array}{ll}
 \mbox{minimize} &  (\Tr V  + w)/2 \\*[1ex]
 \mbox{subject to} &  \left[\begin{array}{cc}
     V &  x \\ x^H & w \end{array}\right] \succeq 0 \\*[1ex]
 & \mbox{$V$ is Toeplitz},
 \end{array}
\EEQ
with variables $w$ and $V\in\herm^n$ (the $n\times n$ Hermitian matrices).
This result can be proved via convex duality and 
semidefinite characterizations of bounded trigonometric
polynomials \cite{dCG12}, or directly by referring to 
Carath\'eodory's  decomposition of positive semidefinite
Toeplitz matrices \cite{TBSR:13}.
More generally, one can consider the atomic norm of the set of matrices
\[
 C = \{ vw^H \in\complex^{n\times m} \mid 
 v \in C_\mathrm{e}, \; \|w\| =1 \}.
\]
The atomic norm for this set, evaluated at a matrix 
$X\in\complex^{n\times m}$, is the optimal value of
the SDP
\BEQ \label{e-Ce-sdp-matrix}
 \begin{array}{ll}
 \mbox{minimize} &  (\Tr V  + \Tr W)/2 \\*[1ex]
 \mbox{subject to} &  \left[\begin{array}{cc}
     V &  X \\ X^H & W \end{array}\right] \succeq 0 \\*[1ex]
 & \mbox{$V$ is Toeplitz},
 \end{array}
\EEQ
with variables $V\in\herm^n$ and $W\in\herm^m$; see \cite{YX14,LC14,Gra15}.
Further extensions, that place restrictions on the parameter 
$\omega$ in the definition~(\ref{e-Ce}), can be found in 
\cite{MCKX14,MCKX15}. 

In this paper we discuss extensions of the SDP 
representations~(\ref{e-Ce-sdp}) and~(\ref{e-Ce-sdp-matrix}) to 
a larger class of atomic norms and gauge functions.
The starting point is the observation that $C_\mathrm e$ can be 
parameterized as
\BEQ \label{e-Ce-def}
 C_\mathrm e = \{ a \mid (\lambda G - F) a = 0, \; \lambda \in
      \mathcal C, \; \|a\|=1\}
\EEQ
where $\mathcal C$ is the unit circle in the complex plane, and
$F$ and $G$ are the $(n-1)\times n$ matrices
\[
 F = \left[\begin{array}{cc} 0 & I_{n-1} \end{array}\right], \qquad
 G = \left[\begin{array}{cc} I_{n-1} & 0 \end{array}\right]. 
\]
We generalize~(\ref{e-Ce-def}) in three ways and derive semidefinite
representations of the corresponding atomic norms.
The first generalization is to replace $\lambda G - F$ with an arbitrary 
matrix pencil.  Second, we allow $\mathcal C$ to be an arbitrary circle or
line in the complex plane, or a segment of a line or a circle.
Third, we replace the normalization $\|a\|=1$ with
a condition of the type $\|Ea\|\leq 1$ where $E$ is not necessarily
full column rank. 
Specific examples of these extensions, with different choices
of $F$, $G$, and $\mathcal C$, are discussed in 
sections~\ref{s-conic-examples}--\ref{s-rational}.

We present direct, constructive proofs, based on elementary matrix 
algebra, of the semidefinite representations of the atomic norms. 
These results are the subject of 
sections~\ref{s-factorization} and~\ref{s-sdp},
and appendix~\ref{s-matrix-fact}.
In section~\ref{s-duality} we derive the convex conjugates of the 
atomic norms and gauge functions, and discuss the relation between 
the dual SDP representations and the Kalman-Yakubovich-Popov 
lemma from linear system theory.
Appendix~\ref{s-slater} contains a discussion of the properties 
of the matrix pencil $\lambda F-G$ that are needed to ensure 
strong duality in the dual problems.
In section~\ref{s-sp-ex} the SDP formulations are illustrated with several
applications in signal processing.

\section{Positive semidefinite matrix factorization}
\label{s-factorization}
Throughout the paper we assume that $F$ and $G$ are complex matrices
of size $p\times n$, and $\Phi$ and $\Psi$ are Hermitian $2\times 2$ 
matrices with $\det\Phi < 0$.  
We define 
\BEQ \label{e-dict}
 \mathcal A = \{ a \in \complex^n \mid
  (\mu G-\nu F)  a = 0, \; (\mu,\nu) \in \mathcal C\},
\EEQ
where 
\BEQ \label{e-mC}
\mathcal C = \left\{(\mu,\nu) \in\complex^2 \mid 
 (\mu,\nu) \neq 0, \; q_\Phi(\mu,\nu) =0, \; 
 q_\Psi(\mu,\nu) \leq 0 \right\}.
\EEQ
Here $q_\Phi$, $q_\Psi$ are the quadratic forms defined by $\Phi$ and
$\Psi$:
\BEQ \label{e-g-def}
q_\Phi(\mu,\nu) 
 = \left[\begin{array}{c} \mu \\ \nu \end{array}\right]^H \Phi 
 \left[\begin{array}{c} \mu \\ \nu \end{array}\right], \qquad
q_\Psi(\mu,\nu) 
 = \left[\begin{array}{c} \mu \\ \nu \end{array}\right]^H \Psi 
 \left[\begin{array}{c} \mu \\ \nu \end{array}\right]. 
\EEQ
The set $\mathcal C$ is a subset of a line or circle in the complex 
plane, expressed in homogeneous coordinates,
as explained in appendix~\ref{s-regions}.

If $\Phi_{11} \neq 0$ or $\Psi_{11} > 0$, then 
$\nu\neq 0$ for all elements $(\mu,\nu) \in \mathcal C$, and we
can simplify the definition of $\mathcal A$ as
\BEQ \label{e-dict-finite}
 \mathcal A = \{ a \in \complex^n \mid
  (\lambda G- F)  a = 0, \; (\lambda,1) \in \mathcal C\}.
\EEQ
If $\Phi_{11} = 0$ and $\Psi_{11} \leq 0$, then the pair $(1,0)$
is also in $\mathcal C$ and the set $\mathcal A$ in~(\ref{e-dict})
is the union of the right-hand side of~(\ref{e-dict-finite})
and the nullspace of $G$. 
Examples of sets $\mathcal A$ are given in 
sections~\ref{s-conic-examples}--\ref{s-rational}.

The purpose of this section is to discuss a semidefinite representation of 
the convex hull of the set of matrices $aa^H$ with $a\in\mathcal A$,
\ie, the set
\BEQ \label{e-KA}
\conv{ \{ aa^H \mid a\in\mathcal A\}} 
= \{ \sum_{k=1}^m a_ka_k^H \mid a\in\mathcal A\}.
\EEQ

\subsection{Conic decomposition}
The key decomposition result (Theorem~\ref{t-decomp}) is known under 
various forms 
in system theory, signal processing, and moment theory 
\cite{KaS:66,KrN:77,GrS:84}.
Our purpose is to give a simple semidefinite formulation that encompasses 
a wide variety of interesting special cases, 
and to present a constructive proof
that can be implemented using the basic decompositions of
numerical linear algebra (specifically, symmetric eigenvalue, 
singular value, and Schur decompositions).

\begin{theorem} \label{t-decomp}
Let $\mathcal A$ be defined by~(\ref{e-dict}) and~(\ref{e-mC}),
where $F$, $G\in\complex^{p\times n}$ and $\Phi$, $\Psi\in\herm^2$ 
with $\det\Phi < 0$.
If $X\in\herm^n$ is a positive semidefinite matrix of rank $r\geq 1$ 
that satisfies 
\BEA
& 
\Phi_{11} FXF^H + \Phi_{21} FXG^H + \Phi_{12} GXF^H + \Phi_{22} GXG^H = 0
& \label{e-phi} \\
& 
\Psi_{11} FXF^H + \Psi_{21} FXG^H + \Psi_{12} GXF^H + \Psi_{22} GXG^H 
\preceq 0,  
& \label{e-psi}
\EEA
then $X$ can be decomposed as 
\BEQ \label{e-X-decomp}
 X = \sum\limits_{k=1}^r a_ka_k^H,
\EEQ
with linearly independent vectors
$a_1$, \ldots, $a_r \in \mathcal A$.
\end{theorem}

\medskip\proof\
We start from any factorization $X=YY^H$ where 
$Y\in\complex^{n\times r}$ has rank $r$.
It follows from Lemma~\ref{l-quad-eq-ineq-general}
in appendix~\ref{s-matrix-fact}, applied to the 
matrices $U=FY$ and $V=GY$, that 
there exist a matrix $W\in\complex^{p\times r}$,
a unitary matrix $Q\in\complex^{r\times r}$,
and two vectors $\mu, \nu \in \complex^r$ such that
\BEQ \label{e-dcmp-pf}
 FYQ = W \diag(\mu), \qquad GYQ = W\diag(\nu), \qquad
(\mu_i, \nu_i) \in \mathcal C, \quad  i=1,\ldots,r.
\EEQ
Choosing $a_k$ equal to the $k$th column of $YQ$
gives the decomposition~(\ref{e-X-decomp}).
\qed\medskip

Viewed geometrically, the theorem says that 
(\ref{e-KA})~is the set of positive semidefinite
matrices $X$ that satisfy~(\ref{e-phi}) and~(\ref{e-psi}).

It is useful to note that the proof of Lemma~\ref{l-quad-eq-ineq-general} 
in the appendix is constructive and gives a simple algorithm, based
on singular value and Schur decompositions, for computing
the matrices $W$, $Q$ and the vectors $\mu$, $\nu$.
In the following three sections we illustrate the decomposition 
in Theorem~\ref{t-decomp} with different choices of 
$F$, $G$, $\Phi$, $\Psi$.

\subsection{Trigonometric polynomials} \label{s-conic-examples}

\paragraph{Complex exponentials}
As a first example, we take $p=n-1$,  
\BEQ \label{e-FG-toep}
F = \left[\begin{array}{cc} 0 & I_{n-1} \end{array}\right], 
\qquad
G = \left[\begin{array}{cc}
 I_{n-1} & 0 \end{array}\right], 
\qquad
 \Phi = \Phi_\mathrm u = 
\left[\begin{array}{cc} 1 & 0 \\ 0 & -1 \end{array} \right], \qquad
 \Psi = 0.
\EEQ
A nonzero pair $(\mu,\nu)$ satisfies 
$q_\Phi(\mu,\nu) = |\mu|^2 - |\nu|^2 =0$  
only if $\mu$ and $\nu$ are nonzero and $\lambda = \mu/\nu$ is on the
unit circle.
The condition $(\lambda G - F)a = 0$ in the definition of $\mathcal A$
gives a recursion
\[
 \lambda a_1 = a_2, \qquad 
 \lambda a_2 = a_3, \qquad \ldots, \qquad \lambda a_{n-1} = a_n.
\]
Defining $\exp(\j\omega)= \lambda$, we find that $\mathcal A$ contains
the vectors
\BEQ \label{e-ak-toep}
 a = c \, (1, e^{\j\omega}, e^{\j 2\omega}, \ldots, e^{\j (n-1)\omega}),
\EEQ
for all $\omega \in [0,2\pi)$ and $c\in\complex$.
The matrix constraints~(\ref{e-phi})--(\ref{e-psi}) reduce to
$FXF^H = GXG^H$,
\ie, $X$ is a Toeplitz matrix.
Theorem~\ref{t-decomp} therefore states that every $n\times n$ positive 
semidefinite Toeplitz matrix can be decomposed as
\BEQ \label{e-carath}
 X = \sum_{k=1}^r |c_k|^2 
 \left[\begin{array}{c} 
 1 \\ e^{\j\omega_k} \\ e^{\j 2\omega_k} \\ \vdots \\ 
  e^{\j (n-1)\omega_k} \end{array}\right]
 \left[\begin{array}{c} 
 1 \\ e^{\j\omega_k} \\ e^{\j 2\omega_k} \\ \vdots \\ 
e^{\j (n-1)\omega_k} \end{array}\right]^H,
\EEQ
with $c_k\neq 0$ and
distinct $\omega_1$, \ldots, $\omega_r$.
This is often called the Carath\'eodory parameterization of
positive semidefinite Toeplitz matrices \cite[page 170]{StM:97}.

For this example, the algorithm outlined in the proof of 
Theorem~\ref{t-decomp} and Lemma~\ref{l-quad-eq-ineq-general}
reduces to the following.
Compute a factorization $X = YY^H$ where $Y\in\complex^{n\times r}$
with rows $y_k^H$, $k=1,\ldots,n$.
Then find a unitary $r\times r$ matrix $\Lambda$ 
that satisfies
\[
\left[\begin{array}{c} y_2^H \\ \vdots \\ y_n^H \end{array}\right]
= \left[\begin{array}{c} y_1^H \\ \vdots \\ y_{n-1}^H \end{array}\right]
 \Lambda,
\]
and compute a Schur decomposition $\Lambda = Q \diag(\lambda)Q^H$.
The eigenvalues give $\lambda_k = \exp(\j\omega_k)$, $k=1,\ldots,r$, and
the columns of $YQ$ are the vectors $a_k$.

\paragraph{Restricted complex exponentials}
Define $F$, $G$, $\Phi$ as in~(\ref{e-FG-toep}), and 
\[
 \Psi = \left[\begin{array}{cc} 0 & -e^{\j\alpha} \\ 
 -e^{-\j\alpha} & 2\cos\beta\end{array}\right]
\]
with $\alpha \in [0, 2\pi)$ and $\beta\in [0,\pi)$.
The elements $a\in\mathcal A$ have the same general 
form~(\ref{e-ak-toep}), with the added constraint that
$\cos\beta \leq \cos(\omega-\alpha)$.
Since we can restrict $\omega$ to the interval 
$[\alpha - \pi, \alpha + \pi]$, this 
is equivalent to $|\omega - \alpha| \leq \beta$.
The constraints~(\ref{e-phi})--(\ref{e-psi}) specify that 
$X$ is Toeplitz and satisfies the matrix inequality
\BEQ \label{e-toep-restricted}
-e^{-\j\alpha} FXG^H  - e^{\j\alpha}GXF^H + 2 (\cos\beta) GXG^H \preceq 0.
\EEQ
The theorem states that a positive semidefinite Toeplitz matrix 
of rank $r$
satisfies~(\ref{e-toep-restricted}) if and only if it can be 
decomposed as~(\ref{e-carath}) with nonzero $c_k$ and
$|\omega_k - \alpha| \leq \beta$ for $k=1,\ldots,r$.

\paragraph{Real trigonometric functions}
Next consider $p=n-1$,
\[
G = \left[\begin{array}{cccccc}
 1 & 0 & 0 & \cdots & 0 & 0  \\
 0 & 2 & 0 & \cdots & 0 & 0  \\
 0 & 0 & 2 & \cdots & 0 & 0  \\
 \vdots & \vdots & \vdots & \ddots & \vdots & \vdots  \\
 0 & 0 & 0 & \cdots & 2 & 0 
 \end{array}\right], \qquad
F = \left[\begin{array}{ccccccc}
 0 & 1 & 0 & \cdots & 0 & 0 & 0  \\
 1 & 0 & 1 & \cdots & 0 & 0 & 0 \\
 0 & 1 & 0 & \cdots & 0 & 0 & 0 \\
 \vdots & \vdots & \vdots &  & \vdots & \vdots &\vdots \\
 0 & 0 & 0 & \cdots & 1 & 0 & 1 
 \end{array}\right], 
\]
and
\[
\Phi = \Phi_\mathrm r = \left[\begin{array}{cc}
 0 & \j \\ -\j & 0 \end{array}\right], \qquad
 \Psi = \Phi_\mathrm u = \left[\begin{array}{cc}
 1 & 0 \\ 0 & -1 \end{array}\right].
\]
A nonzero pair $(\mu,\nu)$ satisfies
$q_\Phi(\mu,\nu) = \j(\bar\mu\nu - \mu\bar \nu) = 0$
and $q_\Psi(\mu,\nu) = |\mu|^2 - |\nu|^2 \leq 0$
only if $\nu \neq 0$  and $\lambda = \mu/\nu$ is real with
$|\lambda| \leq 1$.
The condition $(\lambda G - F)a = 0$ gives a recursion
\[
 \lambda a_1 = a_2, \qquad
 2\lambda a_2 = a_1+a_3, \qquad
 \ldots, \qquad
 2\lambda a_{n-1} = a_{n-2} +a_n. 
\]
If we write $\lambda = \cos\omega$, we recognize the recursion
$2\cos\omega \cos{k\omega} = \cos{(k-1)\omega} + \cos{(k+1)\omega}$
and find that $\mathcal A$ contains the vectors
\[
a = c \, (1,\, \cos\omega,\, \cos{2\omega},\, \ldots, \,
 \cos{(n-1)\omega}),
\]
for all $\omega\in[0,2\pi)$ and all $c$.
With the same $F$ and 
$G = [\begin{array}{cc} 2I_{n-1} & 0 \end{array}]$,
the condition $(\lambda G-F)a = 0$ reduces to
\[
2\lambda a_1 = a_2, \qquad
2\lambda a_2 = a_1 + a_3, \qquad
 \ldots, \qquad
2\lambda a_{n-1} = a_{n-2} + a_n. 
\]
If we write $\lambda = \cos\omega$, the solutions are the vectors
\[
a = c\; (1,  \frac{\sin{2\omega}}{\sin{\omega}}, \,
 \frac{\sin{3\omega}}{\sin{\omega}}, \, \ldots, \, 
 \frac{\sin{n\omega}}{\sin{\omega}}), 
\]
for all $\omega \in [0,2\pi)$ and all $c$.

\paragraph{Trigonometric vector polynomials}
We take $p=(k-1)l$, $n=kl$, 
and replace $F$ and $G$ in~(\ref{e-FG-toep})
with
\[
F = \left[\begin{array}{ccccc} 
0 & I & 0 & \cdots & 0 \\
0 & 0 & I & \cdots & 0 \\
\vdots & \vdots & \vdots & \ddots & \vdots \\
0 & 0 & 0 & \cdots & I \end{array}\right], 
\qquad
G = \left[\begin{array}{ccccc} 
I & 0 & \cdots & 0 & 0\\
0 & I & \cdots & 0 & 0 \\
\vdots & \vdots & \ddots & \vdots & \vdots \\
0 & 0 & \cdots & I & 0 \end{array}\right],
\]
and blocks of size $l\times l$.
Then $\mathcal A$ contains the vectors of the form
\[
a =  (1,\, e^{\j\omega}, \, e^{\j 2\omega}, \, \ldots, \,
 e^{\j  (k-1)\omega} )\otimes c,
\]
for all $c\in\complex^l$ and  $\omega \in [0,2\pi)$,
where $\otimes$ denotes Kronecker product.

\subsection{Polynomials}
\paragraph{Real powers}
Next, define $F$, $G$ as in~(\ref{e-FG-toep}), and
\BEQ \label{e-power-moments}
 \Phi = \Phi_\mathrm r = 
\left[\begin{array}{cc} 0 & \j \\ -\j & 0 \end{array} \right], \qquad
 \Psi = 0.
\EEQ
A pair $(\mu,\nu)$ satisfies $q_\Phi(\mu,\nu) = 0$ if and only if 
$\bar \mu\nu$ is real.
If $(\mu,\nu) \neq 0$, we either have $\nu=0$ and $\mu$ arbitrary,
or $\nu\neq 0$ and $\lambda = \mu/\nu$ real.
The set $\mathcal A$ therefore contains the vectors
\[
a = c \, (1, \lambda, \lambda^2 , \ldots, \lambda^{n-1}), \qquad
a = c \, (0, 0, \ldots, 0, 1)
\]
for all $\lambda\in\reals$ and $c$.
The matrix constraints~(\ref{e-phi})--(\ref{e-psi}) reduce to
$FXG^H = GXF^H$,
\ie, $X$ is a symmetric (real) Hankel matrix.
Hence, a real symmetric positive  semidefinite Hankel matrix of rank
$r$ can be decomposed in one of two forms
\[
 X = \sum_{k=1}^r
 c_k^2
 \left[\begin{array}{c} 1 \\ \lambda_k \\ \vdots \\ \lambda_k^{n-2}
 \\ \lambda_k^{n-1}
 \end{array}\right]
 \left[\begin{array}{c} 1 \\ \lambda_k \\ \vdots \\ \lambda_k^{n-2}
 \\ \lambda_k^{n-1}
 \end{array}\right]^T, \qquad
 X = \sum_{k=1}^{r-1}
 c_k^2
 \left[\begin{array}{c} 1 \\ \lambda_k \\ \vdots \\ \lambda_k^{n-2}
 \\ \lambda_k^{n-1}
 \end{array}\right]
 \left[\begin{array}{c} 1 \\ \lambda_k \\ \vdots \\ \lambda_k^{n-2}
 \\ \lambda_k^{n-1}
 \end{array}\right]^T +  |c_r|^2
 \left[\begin{array}{c} 0 \\ 0 \\ \vdots \\ 0 \\ 1 \end{array}\right]
 \left[\begin{array}{c} 0 \\ 0 \\ \vdots \\ 0 \\ 1 \end{array}\right]^T,
\]
with distinct real $\lambda_k$ and nonzero $c_k$.

\paragraph{Restricted polynomials}
If $\Psi=0$ in~(\ref{e-power-moments}) is replaced by 
\[
 \Psi = \left[\begin{array}{cc}
   2 & -(\alpha + \beta) \\ -(\alpha +\beta) & 2\alpha \beta
  \end{array}\right]
\]
where $-\infty < \alpha < \beta < \infty$, then $\mathcal A$ contains 
all vectors $a = c (1, \lambda, \ldots, \lambda^{n-1})$
with $\lambda \in [\alpha, \beta]$.
The matrix constraints require $X$ to be a real symmetric Hankel matrix
that satisfies
\[
 2FXF^H - (\alpha+\beta) (FXG^H + GXF^H) + 2\alpha\beta GXG^H \preceq 0.
\]

\paragraph{Orthogonal polynomials}
Let $p_0(\lambda)$, $p_1(\lambda)$, $p_2(\lambda)$,~\ldots\ be
a sequence of real polynomials on $\reals$, with $p_i$ of degree $i$.
It is well known that the polynomials 
are orthonormal with respect to an inner product that satisfies the 
property
\BEQ \label{e-shift}
 \langle f(\lambda), \lambda g(\lambda) \rangle = 
 \langle \lambda f(\lambda), g(\lambda)\rangle
\EEQ
(for example, an inner product of the form
$\langle f, g\rangle = 
\int f(\lambda) g(\lambda) w(\lambda) d\lambda$ with $w(\lambda) \geq 0$)
if and only if the polynomials satisfy 
a three-term recursion
\BEQ \label{e-3-term}
 \beta_{i+1}p_{i+1}(\lambda) = 
 (\lambda-\alpha_i) p_i(\lambda) - \beta_i p_{i-1}(\lambda),
\EEQ
where $p_{-1}(\lambda) = 0$ and 
$p_0(\lambda) = 1/d_0$ where $d_0^2 = \langle 1,1\rangle$.
This can be seen as follows \cite{GoK:83}.

Suppose $p_0$, \ldots, $p_{n-1}$ is any set of polynomials, with $p_i$ of
degree $i$.   Then $\lambda p_i(\lambda)$ can be expressed as a linear 
combination of the polynomials $p_0(\lambda)$, \ldots, 
$p_{i+1}(\lambda)$, and therefore
\BEQ \label{e-3term-pf}
\lambda \left[\begin{array}{c}
 p_0(\lambda) \\ p_1(\lambda) \\ \vdots \\ p_{n-2}(\lambda) 
 \end{array}\right]
= \left[\begin{array}{cc} J & \beta_{n-1} e_{n-1} \end{array}\right]
 \left[\begin{array}{c}
 p_0(\lambda) \\ p_1(\lambda) \\ \vdots \\ p_{n-1}(\lambda) 
 \end{array}\right]
\EEQ
for some lower-Hessenberg matrix $J$ 
(\ie, satisfying $J_{ij} = 0$ for $j>i+1$).
Let $\langle \cdot, \cdot \rangle$ be an inner product 
on the space of polynomials of degree $n-1$ or less.
Taking inner products on both sides of~(\ref{e-3term-pf}), we find that
\[ 
 H = JG + \beta_{n-1} e_{n-1}g^T
\]
where  
\[
 H_{ij} = \langle \lambda p_{i-1}(\lambda), p_{j-1}(\lambda)\rangle,
\qquad
 G_{ij} = \langle p_{i-1}(\lambda), p_{j-1}(\lambda)\rangle,
\qquad
 g_j = \langle p_{n-1}(\lambda), p_{j-1}(\lambda)\rangle,
\]
for $i,j = 1, \ldots,n-1$.  The polynomials are orthonormal for the 
inner product if and only if $G=I$ and $g=0$. 
The inner product satisfies the property~(\ref{e-shift})
if and only if $H$ is symmetric.
Hence if the polynomials are orthonormal for an inner product 
that satisfies~(\ref{e-shift}), then $J$ is a symmetric tridiagonal 
matrix.
If we use the notation 
\BEQ \label{e-jacobi}
J = \left[\begin{array}{ccccccc}
  \alpha_0 & \beta_1 & 0 & \cdots & 0 & 0 \\
  \beta_1  & \alpha_1 & \beta_2 & \cdots & 0 & 0 \\
  0        & \beta_2  & \alpha_2 & \cdots & 0 & 0 \\
\vdots & \vdots & \vdots & \ddots & \vdots & \vdots \\
 0 & 0 & 0 & \cdots & \alpha_{n-3} & \beta_{n-2} \\
 0 & 0 & 0 & \cdots & \beta_{n-2} & \alpha_{n-2} 
 \end{array}\right],
\EEQ
the recursion~(\ref{e-3-term}) follows.
Conversely, if the three-term recursion holds, and we define
the inner product by setting $G=I$, $g=0$, then $H$ is symmetric and
the inner product satisfies~(\ref{e-shift}).

Now consider~(\ref{e-dict}) and~(\ref{e-mC}), with $p=n-1$ and
\[
\Phi = \Phi_\mathrm r, \qquad \Psi = 0, \qquad
G = \left[\begin{array}{cc}
I_{n-1} & 0
\end{array}\right], \qquad
F = \left[\begin{array}{cc} 
 J & \beta_{n-1}e_{n-1} \end{array}\right],
\]
where $J$ is the Jacobi matrix~(\ref{e-jacobi}) of a system of 
orthogonal polynomials.
Then $(\mu,\nu) \in \mathcal C$ if and only if either
$\nu\neq 0$ and $\lambda = \mu/\nu \in \reals$, or $\nu=0$.
The set contains the vectors $a$  of the form
\[
a = c \, (p_0(\lambda),\, p_1(\lambda), \, p_2(\lambda),\ldots, 
\, p_{n-1}(\lambda)), \qquad
a = c \, (0,0, \ldots, 0,1)
\]
for all $\lambda \in\reals$. 

\subsection{Rational functions} \label{s-rational}
As a final example, we consider the controllability pencil of a linear 
system:
\BEQ \label{e-state-space}
 G = \left[\begin{array}{cc} I & 0 \end{array}\right], 
 \qquad
 F = \left[\begin{array}{cc} A & B \end{array}\right],
\EEQ
where $A \in \complex^{n_\mathrm s \times n_\mathrm s}$
and $B\in \complex^{n_\mathrm s \times m}$.
With this choice, $\mathcal A$ contains the vectors $a =(x,u)$
that satisfy the equality $(\mu I - \nu A) x = \nu Bu$ for some 
$(\mu, \nu) \in \mathcal C$.
Since $(\mu,\nu) \neq 0$, we either have $\nu= 0$ and $x=0$,
or $\nu\neq 0$ and $((\mu/\nu) I - A) x= Bu$.
If $A$ has no eigenvalues $\lambda$ that satisfy
$(\lambda, 1) \in \mathcal C$,  then $\mathcal A$ contains the vectors 
\[
 a = \left[\begin{array}{c} (\lambda I -A)^{-1}B u \\ u \end{array}\right]
\]
for all $(\lambda, 1)\in\mathcal C$ and all $u\in\complex^m$.
If $\mathcal C$ includes the point $(1,0)$ at infinity, then
$\mathcal A$ also contains the vectors $(0,u)$ for all $u\in\complex^m$.

This can be extended to the controllability pencil of a descriptor system
\[
 G = \left[\begin{array}{cc} E & 0 \end{array}\right], 
 \qquad
 F = \left[\begin{array}{cc} A & B \end{array}\right],
\]
where $E\in\complex^{n_\mathrm s \times n_\mathrm s}$ is possibly
singular.  
With this choice, $\mathcal A$ contains the vectors $a =(x,u)$
that satisfy the equality $(\mu E - \nu A) x = \nu Bu$ for some 
$(\mu, \nu) \in \mathcal C$.
If $\det(\mu E - \nu A) \neq 0$ for all $(\mu,\nu) \in \mathcal C$, 
then $\mathcal A$ contains all vectors
\[
 a = \left[\begin{array}{c} (\lambda E -A)^{-1}B u \\ u \end{array}\right]
\]
for all $(\lambda, 1)\in \mathcal C$ and all $u\in\complex^m$.  
If $(0,1) \in \mathcal C$, then $\mathcal A$ also contains the points 
$(0,u)$ for all $u\in\complex^m$.

\section{Semidefinite representation of gauges and atomic norms} 
\label{s-sdp}

A function $g$ is called a \emph{gauge} if it is convex, 
positively homogeneous ($g(tx) = tg(x)$ for $t>0$), nonnegative, 
and vanishes at the origin
\cite[section 15]{Roc:70}, \cite[chapter 1]{KrN:77}.
Examples are the \emph{(Minkowski) gauges} of nonempty convex sets $C$,  
which are defined as
\[
 g(x) = \inf{\{ t \geq 0 \mid x \in t C\}}.
\]
Conversely, if $g$ is a gauge, then it is the Minkowski gauge of the set
$C = \{ x \mid g(x) \leq 1\}$.
A gauge is a norm if it is defined everywhere, positive except at
the origin, and symmetric ($g(x) = g(-x)$).

The gauge of the convex hull $\conv C$ of a set $C$ can be expressed as
\[
 g(x) = \inf{\{\sum_{k=1}^r \theta_k \mid
  x = \sum_{k=1}^r \theta_k x_k, \; \theta_k \geq 0, \; x_k \in C, \; 
  k=1,\ldots,r \}}.
\]
The minimum is over all possible decompositions of $x$ as a nonnegative
combination of a finite number of elements of $C$.
The gauge of the convex hull of a compact set is 
also called the \emph{atomic norm} associated with the set
\cite{CRPW:12}.  

\subsection{Symmetric matrices} \label{s-gauge}
Let $F$, $G$, $\Phi$, $\Psi$ be defined as in Theorem~\ref{t-decomp}.
We assume that the set $\mathcal C$ defined in~(\ref{e-mC}) is not empty.
In this section we discuss the gauge of the convex hull of the set 
\[
 C = \{ aa^H \in \herm^n \mid a \in \mathcal A, \; \|a\|=1\},
\]
where $\mathcal A$ is defined in~(\ref{e-dict}).
The gauge of the convex hull of $C$ is the function
\BEA
g(X) & = & \inf{\{\sum_{k=1}^r \theta_k \mid 
  X = \sum_{k=1}^r \theta_k a_ka_k^H, \; \theta_k \geq 0, \;
 a_k \in \mathcal A, \; \|a_k\| =1, \; k=1,\ldots, r\}} 
 \label{e-gammaC-def-a} \\
 & = & \inf{\{\sum_{k=1}^r \|a_k\|^2 \mid 
  X = \sum_{k=1}^r a_ka_k^H, \; a_k \in \mathcal A, \; k=1,\ldots, r\}}.
 \label{e-gammaC-def}
\EEA
The second expression follows from the fact that if $a \in \mathcal A$
then $\beta a\in \mathcal A$ for all $\beta$.

The expressions $\sum_k \theta_k$ and $\sum_k \|a_k\|^2$ in these
minimizations take only two possible values: $\Tr X$ if $X$
can be decomposed as in~(\ref{e-gammaC-def-a}) 
and~(\ref{e-gammaC-def}), and $+\infty$ otherwise.
Theorem~\ref{t-decomp} tells us that a decomposition exists
if only if $X$ is positive semidefinite and satisfies the two 
constraints~(\ref{e-phi}),~(\ref{e-psi}).
Therefore
\BEQ \label{e-gammaC-lmi}
g(X) = \left\{ \begin{array}{ll}
  \Tr X & \mbox{$X\succeq 0$, (\ref{e-phi}),~(\ref{e-psi})} \\
 + \infty & \mbox{otherwise.} 
\end{array}\right.
\EEQ

Now consider an optimization problem in which we minimize the sum of a 
function $f:\herm^n \rightarrow\reals$ and the gauge 
defined in~(\ref{e-gammaC-def}) and~(\ref{e-gammaC-lmi}),
\BEQ \label{e-f-plus-gamma}
 \begin{array}{ll}
 \mbox{minimize} & f(X) + g(X).
 \end{array}
\EEQ
If we substitute the definition~(\ref{e-gammaC-def}),
this can be written as 
\BEQ \label{e-f-plus-gauge}
 \begin{array}{ll}
 \mbox{minimize} & f(X) + \sum\limits_{k=1}^r \|a_k\|^2 \\*[1ex]
 \mbox{subject to} & X = \sum\limits_{k=1}^r a_ka_k^H \\*[2ex]
 & a_k \in \mathcal A, \; k=1,\ldots, r.
 \end{array}
\EEQ
The variables are $X$ and the parameters $a_1$, \ldots, $a_r$, and
$r$ of the decomposition of $X$.
This formulation shows that the function $g(X)$ 
in~(\ref{e-f-plus-gamma}) acts as a regularization term that
promotes a structured low rank property in $X$. 
If we substitute the expression~(\ref{e-gammaC-lmi}) we obtain the
equivalent formulation
\BEQ \label{e-gauge-sdp-rep}
\begin{array}{ll}
\mbox{minimize} & f(X) + \Tr X \\
 \mbox{subject to} 
 & \Phi_{11} FXF^H + \Phi_{21} FXG^H + 
        \Phi_{12} GXF^H + \Phi_{22} GXG^H = 0\\
 & \Psi_{11} FXF^H + \Psi_{21} FXG^H + 
     \Psi_{12} GXF^H + \Psi_{22} GXG^H \preceq 0  \\
 & X \succeq 0.
\end{array}
\EEQ
This problem is convex if $f$ is convex. 

A useful generalization of~(\ref{e-gammaC-def}) is the 
gauge of the convex hull of
\[
 C = \{ aa^H \mid a\in \mathcal A, \; \|Ea\| \leq 1\}
\]
where $E$ may have rank less than $n$.
The gauge of $\conv{C}$ is
\BEQ \label{e-gammaCP-def}
g(X) = \inf{\{ \sum_{k=1}^r \theta_k \mid
   X =\sum_{k=1}^r \theta_k a_ka_k^H, \; \theta_k\geq 0, \;
   a_k\in\mathcal A, \; \|Ea_k\| \leq 1, \; k=1,\ldots, r \}}.
\EEQ
The variables $\theta_k$ in this definition can be eliminated
by making the following observation.
Suppose that the directions of the vectors $a_k$ in the decomposition 
of $X$ in~(\ref{e-gammaCP-def}) are given, but not their norms or
the coefficients $\theta_k$.
If $0 < \|Ea_k\| < 1$, we can decrease $\theta_k$ by scaling 
$a_k$ until $\|Ea_k\| = 1$. 
If $Ea_k = 0$, $\theta_k$ can be made arbitrarily small by
scaling $a_k$.
Hence, we obtain the same result if we use $\sqrt \theta_k a_k$ as 
variables and write the infimum as:  
\BEQ \label{e-gammaC-def-weighted}
g(X) 
= \inf{\{ \sum_{k=1}^r \|Ea_k\|^2  \mid X = \sum_{k=1}^r a_ka_k^H, \;
 a_k \in \mathcal A, \; k=1,\ldots, r\}}.
\EEQ
Therefore $g(X) = \sum_k \|Ea_k\|^2 = \Tr(EXE^H)$ if $X$
can be decomposed as in~(\ref{e-gammaC-def-weighted}) and 
$+\infty$ otherwise.
Using Theorem~\ref{t-decomp} we can express this result as
\BEQ \label{e-gammaCP-lmi}
 g(X) 
 = \left\{ \begin{array}{ll}
  \Tr(EXE^H) & \mbox{$X\succeq 0$, (\ref{e-phi}), (\ref{e-psi})} \\
  + \infty & \mbox{otherwise.} \end{array}\right.
\EEQ
Minimizing $f(X) + g(X)$ is equivalent to
the optimization problem
\BEQ \label{e-gauge-sdp-rep-E}
 \begin{array}{ll}
 \mbox{minimize} & f(X) + \sum\limits_{k=1}^r \|Ea_k\|^2  \\*[1ex]
 \mbox{subject to} & X = \sum\limits_{k=1}^r a_ka_k^H \\*[1.5ex]
 & a_k \in \mathcal A, \;\; k=1,\ldots, r,
 \end{array}
\EEQ
with variables $X$  and the parameters $a_1$, \ldots, $a_r$, $r$
of the decomposition of $X$.
When $E^HE=I$ this is the same as~(\ref{e-f-plus-gauge}).
By choosing different $E$ we assign different weights
to the vectors $a_k$. 
Using the expression~(\ref{e-gammaCP-lmi}), the 
problem~(\ref{e-gauge-sdp-rep-E}) can be written as
\BEQ \label{e-gauge-sdp-2}
\begin{array}{ll}
 \mbox{minimize} & f(X) + \Tr{(EXE^H)} \\
 \mbox{subject to} 
 & \Phi_{11} FXF^H + \Phi_{21} FXG^H + 
        \Phi_{12} GXF^H + \Phi_{22} GXG^H = 0\\
 & \Psi_{11} FXF^H + \Psi_{21} FXG^H + 
     \Psi_{12} GXF^H + \Psi_{22} GXG^H \preceq 0  \\
 & X \succeq 0.
\end{array}
\EEQ

\paragraph{Example}
Parametric line spectrum estimation is concerned with fitting
signal models of the form
\BEQ \label{e-sinusoids-in-noise}
 y(t) = \sum_{k=1}^r c_k e^{\j \omega_k t} + v(t),
\EEQ
where $v(t)$ is noise. 
If the phase angles of $c_k$ are independent random variables,
uniformly distributed on $[-\pi,\pi]$, and $v(t)$ is 
circular white noise with $\Expect |v(t)|^2  = \sigma^2$, then
the covariance matrix of $y(t)$ of order $n$ is given by
\BEQ \label{e-line-spec-covariance}
\left[\begin{array}{cccc}
 r_0 & r_{-1} & \cdots & r_{-n+1} \\
 r_1 & r_0   & \cdots & r_{-n+2} \\
 \vdots & \vdots & \ddots & \vdots \\
 r_{n-1} & r_{n-2}   & \cdots & r_0
 \end{array}\right]
=
\sigma^2 I + \sum\limits_{k=1}^r |c_k|^2
 \left[\begin{array}{c} 
  1 \\ e^{\j\omega_k} \\ \vdots \\ e^{\j (n-1)\omega_k}  
 \end{array}\right]
 \left[\begin{array}{c} 
  1 \\ e^{\j\omega_k} \\ \vdots \\ e^{\j (n-1)\omega_k}  
 \end{array}\right]^H,
\EEQ
where $r_k = \Expect{(y(t) \overline{y(t-k)})}$
\cite[section 4.1]{StM:97}\cite[section 12.5]{PrM:96}.
Classical methods, such as MUSIC and ESPRIT, 
are based on the eigenvalue decomposition of an 
estimated covariance matrix.
With the formulation outlined in this section one can solve
related but more general covariance fitting problems, expressed as
\[
 \begin{array}{ll}
 \mbox{minimize} & f(R) + 
 n \sum\limits_{k=1}^r  |c_k|^2 \\
 \mbox{subject to} & 
  R = \sigma^2 I + \sum\limits_{k=1}^r |c_k|^2
 \left[\begin{array}{c} 
  1 \\ e^{\j\omega_k} \\ \vdots \\ 
  e^{\j (n-1)\omega_k}  \end{array}\right]
 \left[\begin{array}{c} 
  1 \\ e^{\j\omega_k} \\ \vdots \\ 
  e^{\j (n-1)\omega_k}  \end{array}\right]^H,
 \end{array}
\]
with variables $R\in\herm^n$, $\sigma^2$, $|c_k|$, $\omega_k$, and $r$,
where $f$ is a convex penalty or indicator function that measures
the quality of the fit between $R$ and the estimated covariance matrix.
This is equivalent to the convex optimization problem
\[
 \begin{array}{ll}
 \mbox{minimize} & f(X + tI) + \Tr X \\ 
 \mbox{subject to} & X \succeq 0, \; t\geq 0 \\
                   & \mbox{$X$ is Toeplitz}. 
 \end{array}
\]
A numerical example is given in section~\ref{s-sp-ex}.

\subsection{Non-symmetric matrices} \label{ss-main}
We define $F$, $G$, $E$, $\Phi$, $\Psi$, and $\mathcal A$ as in the 
previous section, but add the assumption that the matrices $F$, $G$, and
$E$ are block-diagonal:
\BEQ \label{e-GF-diag}
 G  = \left[\begin{array}{cc}
   G_1 & 0 \\ 0 & G_2 \end{array}\right], \qquad
 F  = \left[\begin{array}{cc}
   F_1 & 0 \\ 0 & F_2 \end{array}\right],  \qquad
 E  = \left[\begin{array}{cc}
   E_1 & 0 \\ 0 & E_2 \end{array}\right].  
\EEQ
Here $F_1$, $G_1 \in \complex^{p_1 \times n_1}$ and
$F_2$, $G_2 \in \complex^{p_2 \times n_2}$ 
(possibly with $p_1$ or $p_2$ equal to zero).  
The matrices $E_1$ and $E_2$ have $n_1$ and $n_2$ columns, 
respectively.  
In this section we discuss the function
\[
 h(Y) =  \frac{1}{2} \inf_{V,W} g(\left[\begin{array}{cc}
   V & Y \\ Y^H & W \end{array}\right])
\]
of $Y\in \complex^{n_1\times n_2}$, where $g$ is the function defined 
in~(\ref{e-gammaC-def-weighted}) and~(\ref{e-gammaCP-lmi}).
Using~(\ref{e-gammaC-def-weighted}) we can write $h(Y)$ as
\BEQ
h(Y)
= \inf{ \{ \frac{1}{2} \sum_{k=1}^r (\|E_1v_k\|^2 + \|E_2w_k\|^2) \mid 
  Y = \sum_{k=1}^r v_k w_k^H, \; (v_k,w_k) \in \mathcal A\}},
 \label{e-gauge-nonsymm-d}
\EEQ
while the equivalent characterization~(\ref{e-gammaCP-lmi}) 
shows that $h(Y)$ is the optimal value of the  SDP
\BEQ \label{e-gauge-nonsymm-sdp}
\begin{array}{ll}
\mbox{minimize} & \left(\Tr(E_1VE_1^H) + \Tr(E_2WE_2^H) \right)/2\\*[1ex]
\mbox{subject to} 
 & \Phi_{11} FXF^H + \Phi_{21} FXG^H + 
        \Phi_{12} GXF^H + \Phi_{22} GXG^H = 0\\*[.5ex]
 & \Psi_{11} FXF^H + \Psi_{21} FXG^H + 
     \Psi_{12} GXF^H + \Psi_{22} GXG^H \preceq 0  \\*[.5ex]
 & X = \left[\begin{array}{cc}
    V & Y \\ Y^H & W \end{array}\right] \succeq 0,
\end{array}
\EEQ
with $V$ and $W$ as variables.
This can be seen as an extension of the well-known SDP formulation of 
the trace norm of a rectangular matrix.
If we take $F$ and $G$ to have zero row dimensions (equivalently,
define $\mathcal A = \complex^{n_1} \times \complex^{n_2}$ and
omit the first two constraints in~(\ref{e-gauge-nonsymm-sdp}))
and choose identity matrices for $E_1$ and $E_2$, 
then $h(Y) = \|Y\|_*$, the trace norm of $Y$.

The block-diagonal structure of $F$ and $G$ implies that if
$(v,w)\in\mathcal A$, then $(\alpha  v, \beta w)\in \mathcal A$
for all $\alpha$, $\beta$.
This observation leads to a number of useful equivalent expressions 
for~(\ref{e-gauge-nonsymm-d}).
First, we note that $h(Y)$ can be written as
\BEQ
 h(Y) 
 =  \inf{ \{ \sum_{k=1}^r \|E_1v_k\| \|E_2w_k\| \mid 
  Y = \sum_{k=1}^r v_k w_k^H, \; (v_k,w_k) \in \mathcal A\}}.
 \label{e-gauge-nonsymm-c}  
\EEQ
This follows from the fact 
$\|E_1v_k\|^2 + \|E_2w_k\|^2 \geq 2\|E_1v_k\| \|E_2w_k\|$,
with equality if $\|E_1v_k\| = \|E_2w_k\|$.  
If the decomposition of $Y$ in~(\ref{e-gauge-nonsymm-d})
involves a term $v_k w_k^H$ with $E_1v_k$ and $E_2w_k$ nonzero, 
then replacing $v_k$ and $w_k$ with
\[
 \tilde v_k = \frac{\|E_2w_k\|^{1/2}}{\|E_1v_k\|^{1/2}} v_k, \qquad
 \tilde w_k = \frac{\|E_1v_k\|^{1/2}}{\|E_2w_k\|^{1/2}} w_k
\]
gives another valid decomposition with
\[
 \frac{1}{2} (\|E_1 \tilde v_k\|^2 + \|E_2 \tilde w_k\|^2)
 = \|E_1v_k\| \|E_2 w_k\|
 \leq \frac{1}{2} (\|E_1 v_k\|^2 + \|E_2 w_k\|^2).
\]
If $E_1v_k = 0$  and $E_2w_k \neq 0$, then replacing
$v_k$ and $w_k$ with $\tilde v_k = \alpha v_k$, 
$\tilde w_k = (1/\alpha) w_k$
gives an equivalent decomposition with 
\[
 \frac{1}{2} (\|E_1 \tilde v_k\|^2 + \|E_2 \tilde w_k\|^2)
 = \frac{1}{2\alpha^2}  \|E_2 w_k\|^2 \rightarrow 0
\]
as $\alpha$ goes to infinity.
The same argument applies when $E_1v_k \neq 0$  and $E_2w_k = 0$. 
In all cases, therefore, the two 
expressions~(\ref{e-gauge-nonsymm-d}) and~(\ref{e-gauge-nonsymm-c}) give
the same result.

From~(\ref{e-gauge-nonsymm-c}) we obtain two other useful expressions:
\BEA
h(Y) & = & \inf{ \{ \sum_{k=1}^r \|E_1v_k\| \mid 
  Y = \sum_{k=1}^r v_k w_k^H, \; 
   (v_k, w_k) \in \mathcal A, \; \|E_2w_k\| \leq 1\}}
 \label{e-gauge-nonsymm-a}  \\
& = & \inf{ \{ \sum_{k=1}^r \|E_2w_k\| \mid 
  Y = \sum_{k=1}^r v_k w_k^H, \; 
   (v_k,w_k) \in \mathcal A, \; \|E_1v_k\| \leq 1\}}.
 \label{e-gauge-nonsymm-b}  
\EEA
This again follows from the property that the two components of 
elements $(v_k, w_k)$ in $\mathcal A$ can be scaled independently.
At the optimal decomposition in~(\ref{e-gauge-nonsymm-a}), all
terms in the decomposition satisfy $E_2w_k = 0$ or $\|E_2w_k\|=1$.
In~(\ref{e-gauge-nonsymm-b}), all terms satisfy 
$E_1v_k = 0$ or $\|E_1v_k\|=1$.

A final interpretation of $h$ is
\BEQ
h(Y) = 
\inf{ \{ \sum_{k=1}^r \theta_k \mid 
  Y = \sum_{k=1}^r \theta_k v_k w_k^H, \; \theta_k \geq 0, \; 
   (v_k,w_k) \in \mathcal A, \; \|E_1v_k\| \leq 1,  \; \|E_2w_k\| \leq 1 
\}}.
\label{e-gauge-nonsymm-weighted}
\EEQ
The equivalence with~(\ref{e-gauge-nonsymm-c}) follows from the fact
that if the optimal decomposition of $Y$ 
in~(\ref{e-gauge-nonsymm-weighted}) involves 
the term $v_kw_k^H$, then the norms $\|E_1v_k\|$ and $\|E_2w_k\|$ will be 
either zero or one.
(If $0 < \|E_1v_k\| < 1$ we can decrease $\theta_k$ by scaling $v_k$
until $\|E_1v_k\| = 1$, and similarly for $w_k$.)
The expression~(\ref{e-gauge-nonsymm-weighted}) shows that 
$h(Y)$ is the gauge of the convex hull of the set
\BEQ \label{e-uv-dict2}
\{ vw^H \in\complex^{n_1\times n_2} \mid  
 (v,w) \in \mathcal A, \, \|E_1v\| \leq 1, \; \|E_2w\| \leq 1\}.
\EEQ

The SDP representation of $h$ in~(\ref{e-gauge-nonsymm-sdp})
allows us to reformulate problems 
\BEQ \label{e-f-plus-gamma-nonsymm}
\mbox{minimize} \quad f(Y) + h(Y), 
\EEQ
where $f$ is convex and $h$ is the 
gauge~(\ref{e-gauge-nonsymm-d})--(\ref{e-gauge-nonsymm-weighted}),
as a convex problem.
Minimizing $f(Y) + h(Y)$ is equivalent to
\BEQ \label{e-f+g-nonsymm}
\begin{array}{ll} 
\mbox{minimize} &
 f(Y) +  \sum\limits_{k=1}^r \|E_1v_k\| \|E_2w_k\| \\
\mbox{subject to} &  Y = \sum\limits_{k=1}^r v_k w_k^H \\
& (v_k, w_k) \in \mathcal A, \; k=1,\ldots, r.
\end{array}
\EEQ
Alternatively, one can replace the second term in the objective
with $\sum_k \|E_2w_k\|$ and add constraints $\|E_1v_k \| \leq 1$,
as in
\BEQ \label{e-f+g-nonsymm-2}
\begin{array}{ll} 
\mbox{minimize} &
 f(Y) +  \sum\limits_{k=1}^r \|E_2w_k\| \\
\mbox{subject to} &  Y = \sum\limits_{k=1}^r v_k w_k^H \\[.5ex]
& (v_k, w_k) \in \mathcal A, \; k=1,\ldots, r \\[.5ex]
& \|E_1v_k\| \leq 1, \; k = 1,\ldots, r,
\end{array}
\EEQ
or vice versa.
When $E_1$ and $E_2$ are identity matrices, we can interpret
$h(Y)$ as a convex penalty that promotes a structured low-rank 
property of $Y$.  The outer products $v_kw_k^H$ are constrained
by the set $\mathcal A$; the penalty term in the objective
is the sum of the norms $\|v_kw_k^H\|_2 = \|v_k\| \|w_k\|$.
The matrices $E_1$ and $E_2$ can be chosen to assign a different
weight to different terms $v_k w_k^H$.

Problems~(\ref{e-f+g-nonsymm}) and~(\ref{e-f+g-nonsymm-2}) can be 
reformulated as
\BEQ \label{e-min-f+g-general}
\begin{array}{ll} 
\mbox{minimize} 
& f(Y) +  (\Tr(E_1VE_1^H) + \Tr(E_2WE_2^H))/2 \\[.5ex]
\mbox{subject to} 
& \Phi_{11} FXF^H + \Phi_{21} FXG^H + \Phi_{12} GXF^H + \Phi_{22} GXG^H 
 = 0\\[.5ex]
 & \Psi_{11} FXF^H + \Psi_{21} FXG^H + 
     \Psi_{12} GXF^H + \Psi_{22} GXG^H \preceq 0  \\[.5ex]
 & X = \left[\begin{array}{cc}
    V & Y \\ Y^H & W \end{array}\right] \succeq 0.
\end{array}
\EEQ

\paragraph{Example: column structure}
When $p_2=0$, the matrices $F$ and $G$ in~(\ref{e-GF-diag}) have the
form
$F = [\begin{array}{cc} F_1 & 0 \end{array}]$ and
$G = [\begin{array}{cc} G_1 & 0 \end{array}]$. 
This means that $\mathcal A = \mathcal A_1 \times \complex^{n_2}$ where
\[
 \mathcal A_1
 = \{v\in\complex^{n_1} \mid
    (\mu G_1 - \nu F_1) v = 0, \; (\mu, \nu) \in \mathcal C\}.
\]
There are no restrictions on the $w$-component of elements
$(v,w) \in \mathcal A$.
Problem~(\ref{e-f+g-nonsymm}) simplifies to
\BEQ \label{e-f+g-nonsymm-column}
\begin{array}{ll} 
\mbox{minimize} &
 f(Y) +  \sum\limits_{k=1}^r \|E_1v_k\| \|E_2w_k\| \\
\mbox{subject to} &  Y = \sum\limits_{k=1}^r v_k w_k^H \\
& v_k \in \mathcal A_1, \; k=1,\ldots, r,
\end{array}
\EEQ
and the equivalent semidefinite formulation~(\ref{e-min-f+g-general}) to
\[
\begin{array}{ll} 
\mbox{minimize} 
& f(Y) +  (\Tr(E_1VE_1^H) + \Tr(E_2WE_2^H))/2 \\[.5ex]
\mbox{subject to} 
& \Phi_{11} F_1VF_1^H + \Phi_{21} F_1VG_1^H + 
 \Phi_{12} G_1VF_1^H + \Phi_{22} G_1VG_1^H = 0\\[.5ex]
 & \Psi_{11} F_1VF_1^H + \Psi_{21} F_1VG_1^H + 
     \Psi_{12} G_1VF_1^H + \Psi_{22} G_1VG_1^H \preceq 0  \\[.5ex]
 & \left[\begin{array}{cc}
    V & Y \\ Y^H & W \end{array}\right] \succeq 0.
\end{array}
\]

As an example, we again consider the signal 
model~(\ref{e-sinusoids-in-noise}).
A natural idea  for estimating the parameters $\omega_k$ and $c_k$
is to solve a nonlinear least squares problem
\[
 \mbox{minimize} \quad
 \sum_{t=0}^{n-1} | y_\mathrm m(t) - \sum_{k=1}^r c_k e^{\j \omega_k t}|^2,
\]
where $y_\mathrm m(t)$ is the observed signal.
This problem is not convex and difficult to solve iteratively without
a good starting point \cite[page 148]{StM:97}.
However, suppose that, instead of fixing $r$, we impose a penalty on 
$\sum_k |c_k|$, and consider the optimization problem
\BEQ \label{e-nlls}
\begin{array}{ll}
\mbox{minimize} & \gamma \|y - y_\mathrm m\|^2 + \sum\limits_{k=1}^r |c_k|
  \\[1ex]
\mbox{subject to} 
 & y = \sum\limits_{k=1}^r  c_k
 \left[\begin{array}{c}
  1 \\ e^{\j\omega_k} \\ \vdots \\ e^{\j (n-1)\omega_k} \end{array}
 \right].
\end{array}
\EEQ
The optimization variables are $y$ and the parameters 
$c_k$, $\omega_k$, $r$ in the decomposition of $y$. 
The vector $y_\mathrm m$ has elements $y_\mathrm m(0)$, \ldots, 
$y_\mathrm m(n-1)$.
This is a special case of~(\ref{e-f+g-nonsymm-2}) with
$f(y) = \gamma \|y-y_\mathrm m\|^2$, 
$n_1 = n$, $n_2=1$,
\[
 E_1 = \frac{1}{\sqrt n} I, \qquad E_2 = 1, \qquad
 F_1 = \left[\begin{array}{cc} 0 & I_{n_1-1}\end{array}\right], \qquad
 G_1 = \left[\begin{array}{cc} I_{n_1-1} & 0 \end{array}\right],
\]
and $\Phi = \Phi_\mathrm u$, $\Psi = 0$, so that
$\mathcal A_1$ is the set of all multiples of the vectors 
$(1, e^{\j\omega}, \ldots, e^{\j(n-1)\omega})$.
The problem is therefore equivalent to the convex problem
\[
\begin{array}{ll}
\mbox{minimize} & \gamma \|y-y_\mathrm m\|^2  + (\Tr V)/(2n) + w/2 
 \\*[1ex]
\mbox{subject to} & 
 \left[ \begin{array}{cc} V & y \\ y^H & w \end{array}\right]
 \succeq 0 \\*[2ex]
 & \mbox{$V$ is Toeplitz}.
\end{array}
\]
A related numerical example will be given in section~\ref{s-huber}.

\paragraph{Example: joint column and row structure}
To illustrate the general problem~(\ref{e-f+g-nonsymm}), we 
consider a variation on the previous example.
Suppose we arrange the observations in an $n\times m$ Hankel matrix
\[
 Y_\mathrm m = 
 \left[\begin{array}{ccccc}
  y_\mathrm{m}(0) & y_\mathrm{m} (1)  & \cdots & y_\mathrm{m}(m-1)  \\
  y_\mathrm{m}(1) & y_\mathrm{m} (2)  & \cdots & y_\mathrm{m}(m) \\
 \vdots & \vdots  & & \vdots \\
  y_\mathrm{m}(n-1) & y_\mathrm{m} (n) & \cdots & y_\mathrm{m}(m+n-1)
 \end{array}\right],
\]
and we fit to this matrix a matrix $Y$ with the same Hankel 
structure and with elements $y(t) = \sum_{k=1}^r c_k \exp(\j \omega_k t)$.
We formulate the problem as 
\BEQ \label{e-hankel-structure}
\begin{array}{ll}
\mbox{minimize} & \gamma \|Y - Y_\mathrm m\|_F^2
 + \sum\limits_{k=1}^r |c_k| \\
\mbox{subject to} & 
 Y = \sum\limits_{k=1}^r
 c_k \left[\begin{array}{c}
 1 \\ e^{\j\omega_k} \\ \vdots \\ e^{\j (n-1)\omega_k} \end{array}\right]
 \left[\begin{array}{c}
 1 \\ e^{-\j\omega_k} \\ \vdots \\ e^{-\j (m-1)\omega_k} 
  \end{array}\right]^H.
\end{array}
\EEQ
This is an instance of~(\ref{e-f+g-nonsymm}) with 
$n_1=n$, $n_2 =m$,
$E_1= (1/\sqrt n) I$, 
$E_2= (1/\sqrt m) I$, 
and
\[
G_1 = \left[\begin{array}{cc}
I_{m-1} & 0
\end{array}\right], \qquad
F_1 = \left[\begin{array}{cc}
0 & I_{m-1}
\end{array}\right], \qquad
G_2 = \left[\begin{array}{cc}
0 & I_{n-1}
\end{array}\right], \qquad
F_2 = \left[\begin{array}{cc}
I_{n-1} & 0
\end{array}\right].
\]
With these parameters, the set $\mathcal A$ contains the pairs $(v,w)$
of the form 
\[
 v = \alpha (1, e^{\j\omega}, \ldots, e^{\j (m-1)\omega}), \qquad
 w = \beta (1, e^{-\j\omega}, \ldots, e^{-\j (n-1)\omega}).
\]
The convex formulation is
\[
\begin{array}{ll}
\mbox{minimize} & \gamma \|Y-Y_\mathrm m\|_F^2
 + (\Tr V) / (2n) + (\Tr W)/(2m) \\*[1ex]
\mbox{subject to} & \left[\begin{array}{cc}
    V & Y \\ Y^H & W \end{array}\right] \succeq 0 \\*[1ex]
 &\left[\begin{array}{cc} F_1 & 0 \\ 0 & F_2 \end{array}\right] 
  \left[\begin{array}{cc} V & Y \\ Y^H & W \end{array}\right]
 \left[\begin{array}{cc} F_1 & 0 \\ 0 & F_2 \end{array}\right]^T 
 = \left[\begin{array}{cc} G_1 & 0 \\ 0 & G_2 \end{array}\right] 
  \left[\begin{array}{cc} V & Y \\ Y^H & W \end{array}\right]
 \left[\begin{array}{cc} G_1 & 0 \\ 0 & G_2 \end{array}\right]^T.
\end{array}
\]
An example is discussed in section~\ref{s-huber}.

\section{Duality} \label{s-duality}
In this section we derive the conjugates of the gauge functions defined in 
section~\ref{s-sdp} and show that they can be interpreted as 
indicator functions of sets of nonnegative or bounded generalized
polynomials.  This gives a useful interpretation of the dual problems 
for~(\ref{e-f-plus-gamma}) and~(\ref{e-f-plus-gamma-nonsymm}).

We assume that the subset of the complex
plane represented by $\mathcal C$ in~(\ref{e-mC}) is one-dimensional, 
\ie, $\mathcal C$ is not a singleton and not the empty set.  
Equivalently, the inequality $q_\Psi(\mu,\nu) \leq 0$
in the definition is either redundant (and $\mathcal C$ represents a line
or circle), or it is not redundant and then there exist elements of 
$\mathcal C$ with $q_\Psi(\mu, \nu) < 0$. 
When stating and analyzing the dual problems, we will need to
distinguish these two cases ($q_\Psi(\mu,\nu) \leq 0$
is redundant or not).  
For the sake of brevity we only give the formulas for the case where 
the inequality is not redundant.
The dual problems for the other case follow by setting $\Psi=0$
and making obvious simplifications.

We also assume that $\mu G - \nu F$ has full row rank 
($\Rank(\mu G - \nu F) = p$) for all nonzero $(\mu,\nu)$).  
This condition will serve as a 
`constraint qualification' that guarantees strong duality.

\subsection{Symmetric matrix gauge} \label{s-duality-symm}
We first consider the conjugate of the function $g$ defined 
in~(\ref{e-gammaCP-lmi}).  The conjugate is defined as 
\[
 g^*(Z) = \sup_{X} {(\Tr(XZ) - g(X))},
\]
\ie, the optimal value of the SDP
\BEQ \label{e-g-sdp}
\begin{array}{ll}
\mbox{maximize} &  \Tr{((Z-E^HE)X)} \\
\mbox{subject to} 
& X\succeq 0 \\
& \Phi_{11} FXF^H + \Phi_{21} FXG^H + \Phi_{12} GXF^H + \Phi_{22} GXG^H 
 = 0 \\
& \Psi_{11} FXF^H + \Psi_{21} FXG^H + \Psi_{12} GXF^H + \Psi_{22} GXG^H 
\preceq 0. 
\end{array}
\EEQ
The dual of this problem is
\BEQ \label{e-g-conj}
\begin{array}{ll}
\mbox{minimize} & 0 \\
\mbox{subject to}  &
 Z - \left[\begin{array}{c} F \\ G \end{array}\right]^H
  (\Phi \otimes P + \Psi \otimes Q) 
 \left[\begin{array}{c} F \\ G \end{array}\right] \preceq E^HE \\
 & Q \succeq 0,
\end{array}
\EEQ
with variables $P, Q\in \herm^p$.
It is shown in appendix~\ref{s-slater} that strong duality holds 
(under the assumptions listed at the top of section~\ref{s-duality}).

If strong duality holds, then $g^*(Z)$ is the optimal value 
of~(\ref{e-g-conj}), \ie, equal to zero if there exist $P$, $Q$
that satisfy the constraints in~(\ref{e-g-conj}), and $+\infty$
otherwise.  We now show that this can be expressed as 
\BEQ \label{e-g-conj-pos-pol}
 g^*(Z) = \left\{\begin{array}{ll}
   0 & a^H Z a \leq \|Ea\|^2 \mbox{\ for all $a\in \mathcal A$} \\
   +\infty  & \mbox{otherwise.}
\end{array}\right.
\EEQ
Suppose $P$ and $Q$ are feasible in~(\ref{e-g-conj}).
Consider any $a\in\mathcal A$ and $(\mu,\nu) \in \mathcal C$ 
with $\mu Ga = \nu Fa$. Define $y = (1/\nu)Ga$ if
$\nu \neq 0$ and $y = (1/\mu)Fa$ otherwise.
Then 
\BEAS
a^H Za - \|Ea\|^2 
& \leq & \left[\begin{array}{c} Fa \\ Ga \end{array}\right]^H
  (\Phi \otimes P + \Psi \otimes Q) 
  \left[\begin{array}{c} Fa \\ Ga \end{array}\right] \\
& = & \left[\begin{array}{c} \mu y \\ \nu y \end{array}\right]^H
  (\Phi \otimes P + \Psi \otimes Q) 
\left[\begin{array}{c} \mu y \\ \nu y \end{array}\right] \\
& = & (y^H P y) q_\Phi(\mu,\nu) + (y^HQy) q_\Psi(\mu,\nu)  \\
& \leq & 0.
\EEAS
The last line follows from $Q\succeq 0$ and $q_\Phi(\mu,\nu) = 0$,
$q_\Psi(\mu,\nu) \leq 0$.
Conversely, if problem~(\ref{e-g-conj}) is infeasible,
then the optimal value is $+\infty$ and, since strong duality holds,
there exist matrices $X$ that are feasible 
for~(\ref{e-g-sdp}) with $\Tr((Z-E^HE)X) > 0$.
Applying Theorem~\ref{t-decomp}  we see that there exist
$a_1, \ldots, a_r \in\mathcal A$ with 
\[
\sum_{k=1}^r (a_k^HZa_k - \|Ea_k\|^2) > 0.
\]
Therefore $a_k^H Z a_k > \|Ea_k\|^2$ for at least one $a_k$.

The interpretation of the conjugate gives useful insight in 
problem~(\ref{e-f-plus-gamma}), where $g$ is
defined in~(\ref{e-gammaCP-lmi}). 
The dual problem is
\[
 \mbox{maximize} \quad -f^*(Z) - g^*(-Z). 
\]
Expanding $g^*(-Z)$ using~(\ref{e-g-conj}) gives  the equivalent problem
\BEQ\label{e-f-plus-gauge-daul-kyp}
\begin{array}{ll}
\mbox{maximize} & -f^*(Z) \\
\mbox{subject to} & 
 -Z - \left[\begin{array}{c} F \\ G \end{array}\right]^H
  (\Phi \otimes P + \Psi \otimes Q) 
 \left[\begin{array}{c} F \\ G \end{array}\right] \preceq E^HE \\
 & Q \succeq 0,
\end{array}
\EEQ
with variables $Z$, $P$, $Q$,  
and using the expression~(\ref{e-g-conj-pos-pol}) we can put the 
constraints in this problem more succinctly as 
\BEQ \label{e-f-plus-gauge-dual}
\begin{array}{ll}
\mbox{maximize} & -f^*(Z) \\
\mbox{subject to} & \|Ea\|^2 + a^HZ a \geq 0 
 \quad\mbox{for all $a\in\mathcal A$}.
\end{array}
\EEQ
This last form leads to an interesting set of optimality conditions.
Suppose $X$ and $Z$ are feasible for~(\ref{e-gauge-sdp-rep-E})
and~(\ref{e-f-plus-gauge-dual}), respectively. 
Then
\BEAS
 f(X) + \sum_{k=1}^r  \|Ea_k\|^2
 & \geq &  -f^*(Z) + \Tr(XZ) + \sum_{k=1}^r \|Ea_k\|^2\\
 & = &  -f^*(Z) + \sum_{k=1}^r (\|Ea_k\|^2 + a_k^H Za_k) \\
 & \geq &  -f^*(Z).
\EEAS
The first inequality follows by definition of $f^*(Z) = 
\sup_{X}{(\Tr(ZX) - f(X))}$, and the second and third line from
primal and dual feasibility.
If $X$ and $Z$ are optimal and strong duality holds, then $f(X) + \sum_k \|Ea_k\|^2 = -f^*(Z)$.
This is only possible if $f(X) + f^*(Z) = \Tr(XZ)$
and
\[
  \|Ea_k\|^2 +a_k^H Z a_k = 0, \quad k=1,\ldots, r.
\]
Hence only the vectors $a\in\mathcal A$ at which the inequality 
in~(\ref{e-f-plus-gauge-dual}) is active, can be used to form an
optimal $X = \sum_k a_ka_k^H$.

\paragraph{Example: Generalized Kalman-Yakubovich-Popov lemma}
When specialized to the controllability pencil~(\ref{e-state-space}), 
the equivalence between the constraints 
in~(\ref{e-f-plus-gauge-dual}) 
and~(\ref{e-f-plus-gauge-daul-kyp}) is known as the
(generalized) Kalman-Yakubovich-Popov lemma 
\cite{Kal:63,Yak:62,Pop:62,Sch:06,IwH:05}.

We assume that $A$ has no eigenvalues $\lambda$ with $(\lambda, 1)
\in \mathcal C$, and that the pair $(A,B)$ is controllable, so 
the pencil satisfies the rank condition
that $\Rank(\lambda F - G) = n_\mathrm s$ for all $\lambda$.
The dual problem~(\ref{e-f-plus-gauge-dual})
becomes
\[
\begin{array}{ll}
\mbox{maximize} & -f^*(Z) \\
\mbox{subject to} & 
\mathcal F(\lambda, Z) \succeq 0 \quad 
 \mbox{for all $(\lambda, 1)\in \mathcal C$} \\
 & M_{22} + Z_{22} \succeq 0 \quad \mbox{if $(1,0) \in \mathcal C$}
\end{array}
\]
where 
\[
\mathcal F(\lambda, Z) = 
\left[\begin{array}{c} (\lambda I - A)^{-1} B \\ I \end{array}\right]^H
\left[\begin{array}{cc}  M_{11} + Z_{11} & M_{12} + Z_{12} 
  \\ M_{21} +Z_{21} & M_{22} + Z_{22} 
 \end{array}\right]  
\left[\begin{array}{c} (\lambda I - A)^{-1} B \\ I \end{array}\right]
\]
and $M= E^HE$.
The function $\mathcal F$  is called the 
\emph{Popov function} with central matrix $M+Z$ \cite{IOW:99,HSK:99}.

\subsection{Non-symmetric matrix gauge}\label{s-duality-nonsymm}
Next we consider the conjugate of the 
gauge defined in~(\ref{e-gauge-nonsymm-d})--(\ref{e-gauge-nonsymm-a}).
We have
\[
 h^*(Z) = \sup_Y{(\Tr(Z^TY) - h(Y))}
\]
where $h(Y)$ is the  optimal value of~(\ref{e-gauge-nonsymm-sdp}).
Therefore $h^*(Z)$ is the optimal value of the SDP
\BEQ \label{e-g-conj-nonsymm-prim}
\begin{array}{ll} 
\mbox{maximize} & \displaystyle \frac{1}{2} \Tr(\left[\begin{array}{cc}
  -E_1^H E_1 & Z \\ Z^H & -E_2^HE_2 \end{array}\right] X) \\*[2ex]
\mbox{subject to} &
\Phi_{11} FXF^H + \Phi_{21} FXG^H + \Phi_{12}GXF^H + \Phi_{22} GXG^H = 0
\\
& \Psi_{11} FXF^H + \Psi_{21} FXG^H + \Psi_{12}GXF^H + \Psi_{22} GXG^H 
\preceq 0 \\
& X \succeq 0.
\end{array}
\EEQ
The dual of this problem is
\BEQ \label{e-g-conj-nonsymm}
\begin{array}{ll}
\mbox{minimize} & 0 \\
\mbox{subject to} & 
\left[\begin{array}{cc}
 0 & Z \\ Z^H & 0 \end{array}\right]
  - \left[\begin{array}{c} F \\ G \end{array}\right]^H
  (\Phi \otimes P + \Psi \otimes Q) 
 \left[\begin{array}{c} F \\ G \end{array}\right] 
 \preceq \left[\begin{array}{cc} 
 E_1^HE_1 & 0 \\ 0 & E_2^H E_2 \end{array}\right] \\
& Q \succeq 0.
\end{array}
\EEQ
As in the previous section, it follows from appendix~\ref{s-slater}
that strong duality holds.  Therefore $h^*(Z)$ 
is equal to the optimal value of~(\ref{e-g-conj-nonsymm}), \ie, zero
if there exists $P$ and $Q$ that satisfy the constraints of this problem,
and $+\infty$ otherwise.
This will now be shown to be equivalent to
\BEA
 h^*(Z) & = & \left\{\begin{array}{ll}
   0 & \Re{(v^H Z w)} \leq (\|E_1v\|^2 + \|E_2w\|^2)/2 
    \quad \mbox{for all $(v,w)\in\mathcal A$} \\
  +\infty & \mbox{otherwise} \end{array} \right. \nonumber \\
 & = & \left\{\begin{array}{ll}
   0 & \Re{(v^H Z w)} \leq \|E_1v\| \|E_2w\| 
    \quad \mbox{for all $(v,w)\in\mathcal A$} \\
  +\infty & \mbox{otherwise.} \end{array} \right.  
 \label{e-g-conj-nonsymm-pos}
\EEA
To see this, first assume $P$ and $Q$ are feasible 
in~(\ref{e-g-conj-nonsymm}), and $a = (v,w) \in \mathcal A$ 
satisfies $(\mu G - \nu F)a=0$ with $(\mu,\nu)\in\mathcal C$. Then
\BEAS
 v^H Z w + w^H Z^H v - \|E_1u\|^2 - \|E_2v\|^2  
& \leq & \left[\begin{array}{c} Fa \\ Ga \end{array}\right]^H
 (\Phi \otimes P + \Psi \otimes Q) 
 \left[\begin{array}{c} Fa \\ Ga \end{array}\right] \\ 
& = & (y^H Py) q_\Phi(\mu,\nu) + (y^HQy) q_\Psi(\mu,\nu) \\
& \leq & 0,
\EEAS
where we defined $y = (1/\nu) Ga$ if $\nu \neq 0$ and $y=(1/\mu)Fa$ 
otherwise.
Conversely, if problem~(\ref{e-g-conj-nonsymm}) is infeasible,
then~(\ref{e-g-conj-nonsymm-prim}) is unbounded above,
so there exists a feasible $X$ with positive objective value.
If we decompose $X$ as in Theorem~\ref{t-decomp}, 
with $a_k = (v_k,w_k)$, we find that
\BEAS
0 & < &
 \Tr( \left[\begin{array}{cc} 
  -E_1^H E_1 & Z \\ Z^H & -E_2^HE_2 \end{array}\right]
\sum_{k=1}^r \left[\begin{array}{c} v_k \\ w_k \end{array}\right]
\left[\begin{array}{c} v_k \\ w_k \end{array}\right]^H) \\
& = & 
 \sum_{k=1}^r (v_k^HZw_k + w_k^H Z^H v_k - \|E_1 v_k\|^2 - \|E_2w_k\|^2)
\EEAS
so at least one term in the sum is positive.  
The second expression for $h^*(Z)$ in~(\ref{e-g-conj-nonsymm-pos})
follows from the block diagonal structure of $F$ and $G$. 

The interpretation of the conjugate $h^*$ can be applied to interpret
the dual of~(\ref{e-f-plus-gamma-nonsymm}), \ie, 
\[
\begin{array}{ll}
\mbox{maximize} & -f^*(Z) - h^*(-Z).
\end{array}
\]
Substituting the expression~(\ref{e-g-conj-nonsymm}) for $h^*(-Z)$, 
one can write this as
\[
\begin{array}{ll}
\mbox{maximize} & -f^*(Z) \\
\mbox{subject to} & 
 \left[\begin{array}{cc}
  0 & -Z \\ -Z^H  & 0 \end{array}\right]
  - \left[\begin{array}{c} F \\ G \end{array}\right]^H
  (\Phi \otimes P + \Psi \otimes Q) 
 \left[\begin{array}{c} F \\ G \end{array}\right] 
 \preceq \left[\begin{array}{cc} E_1^H E_1 & 0 \\ 0 & E_2^H E_2
 \end{array}\right] \\
 & Q \succeq 0,
\end{array}
\]
with variables $Z$, $P$, $Q$.
Substituting~(\ref{e-g-conj-nonsymm-pos}) we obtain
\[
\begin{array}{ll}
\mbox{maximize} & -f^*(Z) \\
\mbox{subject to} & 
 \Re{(v^H Z w)} \leq \|E_1v\| \|E_2w\| 
 \quad \mbox{for all $(v,w)\in\mathcal A$}.
\end{array}
\]
As in the previous section, the primal-dual optimality conditions
provide a useful set of complementary slackness relations
between primal optimal $Y$ and dual optimal $Z$.
The optimal $Y$ can be decomposed as $Y = \sum_k v_k w_k^H$ with
elements $(v_k,w_k) \in \mathcal A$ at which 
$\Re{(v_k^H Z w_k)} = \|E_1 v_k\| \|E_2 w_k\|$.

\paragraph{Example}
Suppose
$A \in \complex^{n_\mathrm s\times n_\mathrm s}$, 
$B \in \complex^{n_\mathrm s\times m}$, 
$C\in\complex^{l\times n_\mathrm s}$, 
$D\in\complex^{l\times m}$ 
are matrices in a state-space model,
and $A$ has no eigenvalues that satisfy $(\lambda, 1) \in \mathcal C$.
We take $p_1=0$, $n_1 = l$, $p_2 = n_\mathrm s$, 
$n_2 = n_\mathrm s + m$,
\[
 G_2 = \left[\begin{array}{cc} I & 0 \end{array}\right], \qquad
 F_2 = \left[\begin{array}{cc} A & B \end{array}\right], \qquad
 E_1= I, \qquad 
 E_2 = \left[\begin{array}{cc} 0 & I \end{array}\right].
\]
With this choice of parameters, $\mathcal A = \complex^l \times
\mathcal A_2$, where $\mathcal A_2$ contains the vectors of the form
\[
 w = \left[\begin{array}{c} 
 (\lambda I - A)^{-1}B u \\ u \end{array}\right]
\]
for all $u\in\complex^{m}$ and all $(\lambda,1)\in\mathcal C$,
plus the vectors $(0,u)$ if $(0,1) \in\mathcal C$.
Since $v$ is arbitrary and $E_1 = I$,  the inequality
in~(\ref{e-g-conj-nonsymm-pos}) reduces to
$\|Zw\|_2 \leq \|E_2w\|$ for all $w\in\mathcal A_2$.
If $Z$ is partitioned as $Z = [\begin{array}{cc}  C & D\end{array}]$,
this is equivalent to a bound on the transfer function 
\BEQ \label{e-tf-bound}
 \|D + C(\lambda I - A)^{-1}B \|_2 \leq 1 \quad \mbox{
 for all $(\lambda, 1) \in \mathcal C$}, \qquad
 \|D\|_2 \leq 1 \quad \mbox{if $(1,0) \in \mathcal C$}.
\EEQ

\section{Examples}
\label{s-sp-ex}
The formulations in section~\ref{s-sdp} will now be illustrated
with a few examples from signal processing.
The convex optimization problems in the examples were solved 
with CVX \cite{GrB:14}.

\subsection{Line spectrum estimation by Toeplitz covariance fitting}
\label{s-ex-covariance}
\begin{figure}
\centering  
\includegraphics[width=.5\linewidth]{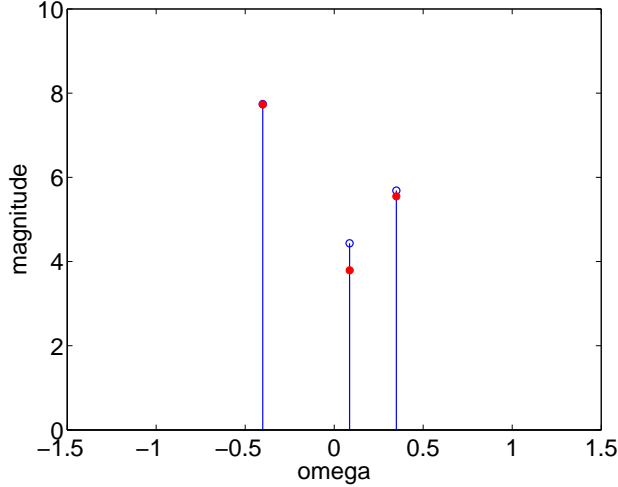}
\caption{Line spectrum estimation by Toeplitz covariance fitting 
(section~\ref{s-ex-covariance}).
The red dots represent the frequencies and magnitudes of the true model.
The blue lines show the estimated parameters obtained by
solving~(\ref{e-cov-fitting}).} \label{f-cov-fitting}
\end{figure}

In this example we fit a covariance matrix of the 
form~(\ref{e-line-spec-covariance}) to an estimated covariance 
matrix~$R_\mathrm m$. 
The estimate $R_\mathrm m$ is constructed from $N=150$ samples 
of the time series $y(t)$ defined in~(\ref{e-sinusoids-in-noise}),
with $r=3$, and frequencies $\omega_k$ and magnitudes $|c_k|$ shown in
figure~\ref{f-cov-fitting}. The noise is Gaussian white noise with variance
$\sigma^2 = 64$.  
The sample covariance matrix is constructed as
\[
R_\mathrm m = \frac{1}{N-n+1} Y Y^H
\]
where $Y$ is the $n\times (N-n+1)$ Hankel matrix 
with $y(1)$, \ldots, $y(N-n+1)$ in its first row.
To estimate the model parameters we solve the optimization problem
\BEQ \label{e-cov-fitting}
 \begin{array}{ll}
 \mbox{minimize} & \gamma \|R-R_\mathrm{m}\|_2 + 
   \sum\limits_{k=1}^r  |c_k|^2 \\
 \mbox{subject to} & R = \sigma^2 I + \sum\limits_{k=1}^r |c_k|^2
 \left[\begin{array}{c} 
  1 \\ e^{\j\omega_k} \\ \vdots \\ 
  e^{\j (n-1)\omega_k}  \end{array}\right]
 \left[\begin{array}{c} 
  1 \\ e^{\j\omega_k} \\ \vdots \\ 
  e^{\j (n-1)\omega_k}  \end{array}\right]^H, %\\
 \end{array}
\EEQ
with variables $|c_k|^2$, $\omega_k$, $r$, and $R$.
The norm $\|\cdot\|_2$ in the objective is the spectral norm.
The regularization parameter $\gamma$ is set to $0.25$.
This problem is equivalent to the convex problem 
\[
 \begin{array}{ll}
 \mbox{minimize} & \gamma \|tI + X -R_\mathrm{m}\|_2 
     + (1/n) \Tr X\\
 \mbox{subject to} & X\succeq 0, \quad t\geq 0 \\
    &  FXF^H - GXG^H = 0 
 \end{array}
\] 
with variables $X$ and $t$, and $F$ and $G$ defined in~(\ref{e-FG-toep}).
As can be seen from Figure~\ref{f-cov-fitting},
the recovered parameters $\omega_k$ and $|c_k|$ are quite accurate,
despite the very low signal-to-noise ratio.
The estimated noise variance $t$ is $79.6$.

The semidefinite optimization approach allows us to fit
a covariance matrix with the structure prescribed 
in~(\ref{e-line-spec-covariance}) to a sample covariance 
matrix that may not be Toeplitz or positive semidefinite.
The formulation can also be extended to applications where the noise
$v(t)$ is modeled as a moving-average process, by combining it with
the formulation in \cite{Geo:06}.

\subsection{Line spectrum estimation by penalty approximation}
\label{s-huber}
\begin{figure}
\centering
\includegraphics[width=.6\linewidth]{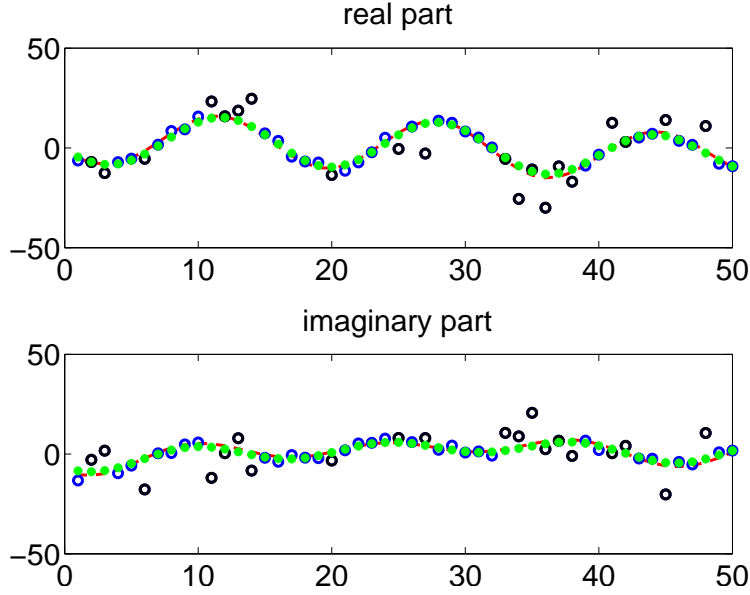} 
\caption{The input data for the example in section~\ref{s-huber}.
The red dashed lines show the exact, noise-free signal.
The blue and black circles show the exact signal corrupted by Gaussian 
white noise, plus a few larger errors in 20 positions.
The green circles show the reconstructed signal.}
\label{f-huber-data}
\end{figure}

\begin{figure}
\centering
\includegraphics[width=.49\linewidth]{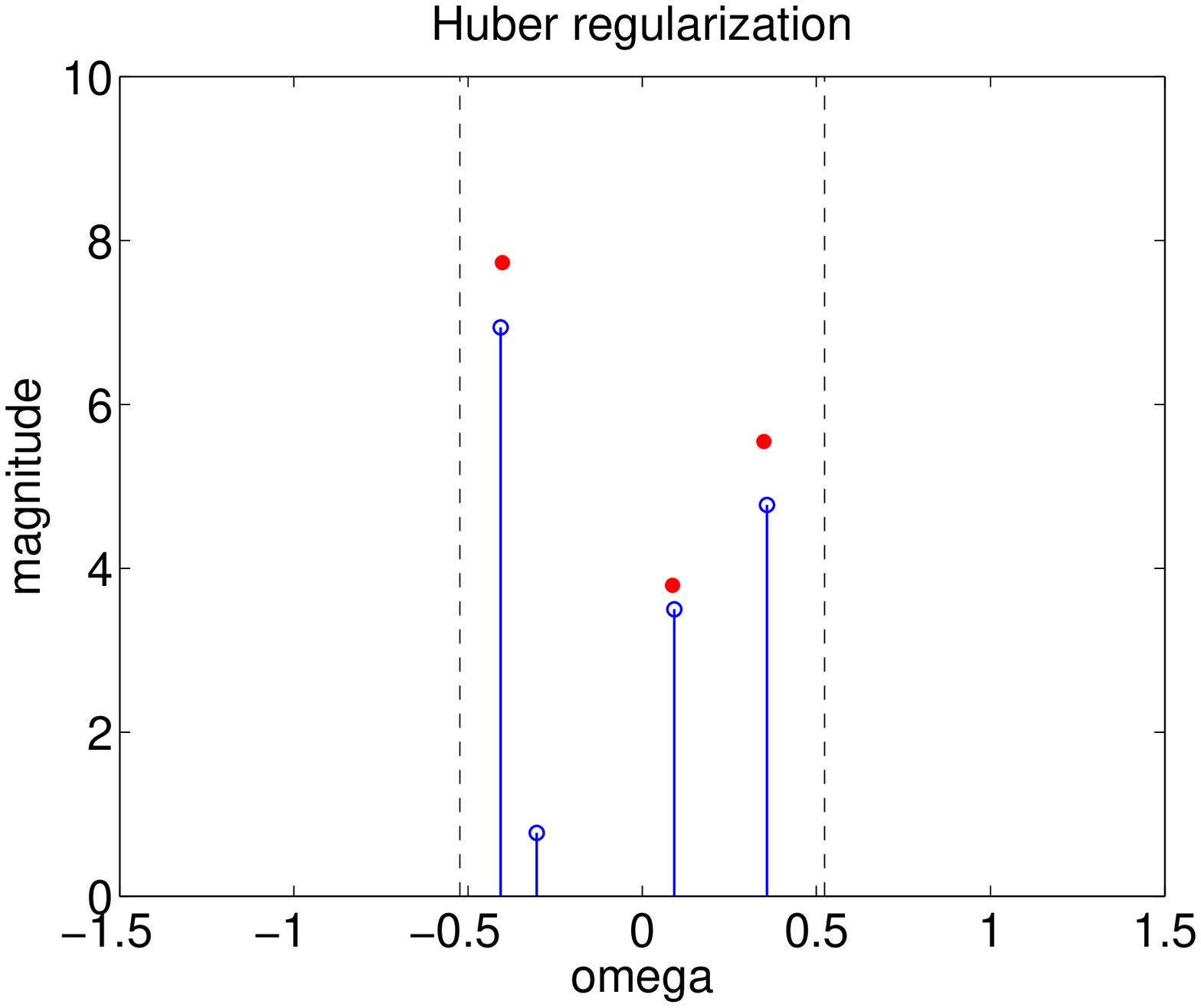}
\hfill
\includegraphics[width=.49\linewidth]{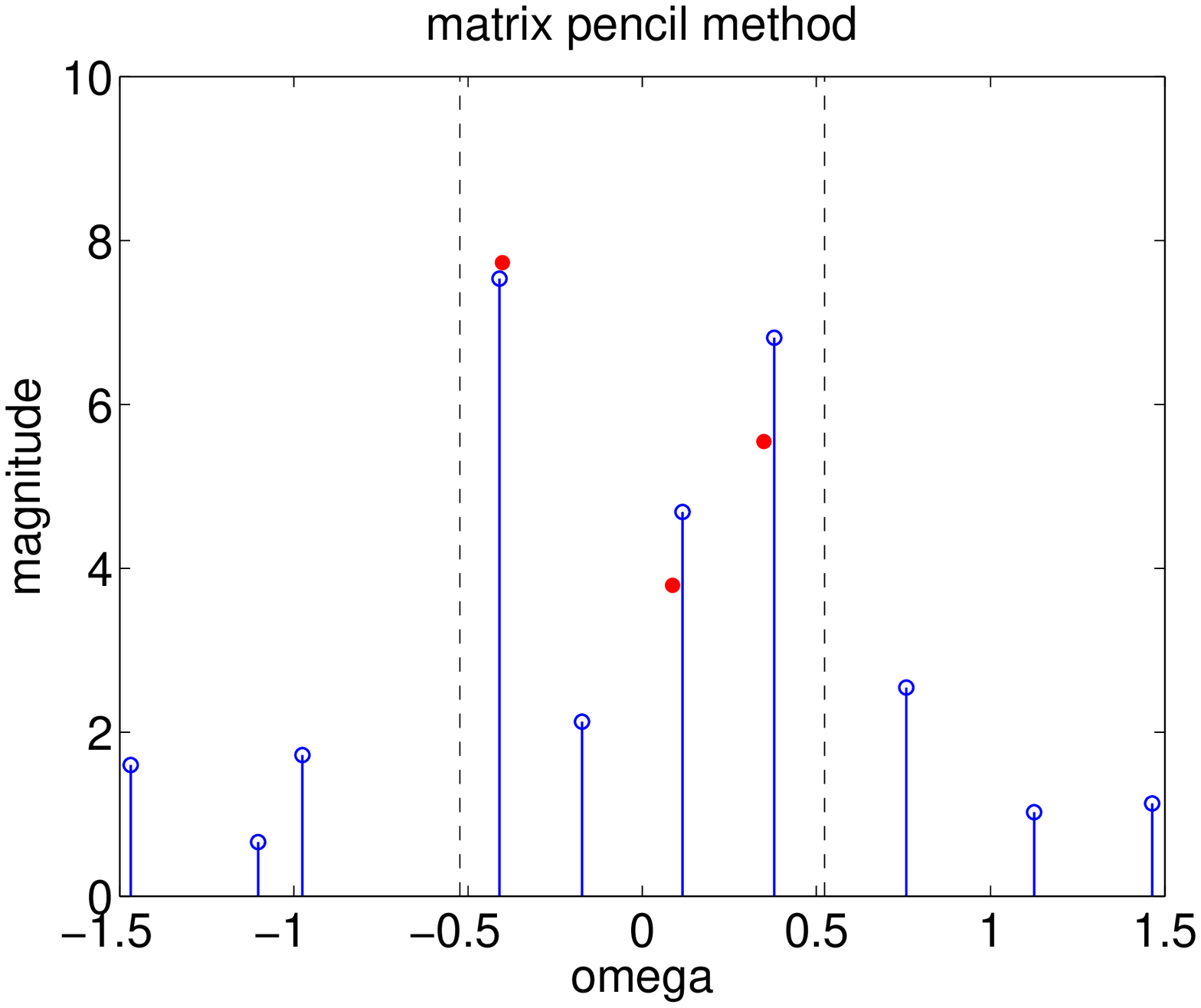}
\caption{Line spectrum models estimated from the signal in 
Figure~\ref{f-huber-data} by solving the optimization 
problem~(\ref{e-huber}) (left) and 
the matrix pencil method (right).}
\label{f-huber}
\end{figure}

This example is a variation on problem~(\ref{e-nlls}).
We take $n=50$ consecutive measurements of the signal
defined in~\eqref{e-sinusoids-in-noise}.  There are 
three sinusoids with frequencies and magnitudes shown in 
figure~\ref{f-huber}.
The noise $v(t)$ is a superposition of white noise and 
a sparse corruption of $20$ elements  (see Figure~\ref{f-huber-data}).
The model parameters are estimated by solving the problem
\BEQ \label{e-huber}
\begin{array}{ll}
\mbox{minimize} & 
\gamma \sum\limits_{i=1}^n \phi(y_i - y_{\mathrm m,i}) + 
  \sum\limits_{k=1}^r |c_k| \\
\mbox{subject to} 
 & y = \sum\limits_{k=1}^r  c_k
 \left[\begin{array}{c}
  1 \\ e^{\j\omega_k} \\ \vdots \\ e^{\j (n-1)\omega_k} \end{array}
 \right] \\
& |\omega_k| \leq \omega_\mathrm c, \quad k = 1,\ldots,r,
\end{array}
\EEQ
where $\phi$ is the Huber penalty, 
$\gamma = 0.071$, 
and $\omega_\mathrm c = \pi/6$.
The variables in this problem are the $n$-vector $y$, and
the parameters $r$, $c_k$, $\omega_k$ in the decomposition of $y$.
The problem is equivalent to the convex problem
\BEQ \label{e-huber-sdp}
\begin{array}{ll}
\mbox{minimize} & 
\gamma \sum\limits_{i=1}^n \phi(y_i-y_{\mathrm m,i})  + 
 (\Tr V)/(2n) + w/2 \\*[1ex]
\mbox{subject to} & 
 \left[ \begin{array}{cc} V & y \\ y^H & w \end{array}\right]
 \succeq 0 \\*[2ex]
 &  FVF^H - GVG^H = 0 \\
 & -FVG^H - GXF^H + 2 (\cos \omega_\mathrm c) GVG^H \preceq 0
\end{array}
\EEQ
with $F$ and $G$ defined in~(\ref{e-FG-toep}).
The variables are the $n$-vector $y$,
the Hermitian $n\times n$ matrix $V$, and the scalar $w$.
The results are shown in Figure~\ref{f-huber}.
The second figure shows the result of a simple 
implementation of the matrix pencil method 
with a $30\times 21$ Hankel matrix constructed from the 
measurements~\cite{HS88}.  
The comparison shows the importance of the prior frequency constraints
in the formulation~(\ref{e-huber}).

It is interesting to note that 
problem~(\ref{e-huber}) can be equivalently formulated as 
\BEQ \label{e-huber-matrix}
\begin{array}{ll}
\mbox{minimize} & 
\gamma \sum\limits_{i=1}^{n} \phi(y_i - y_{\mathrm m,i})
 + \sum\limits_{k=1}^r |c_k| \\
\mbox{subject to} & 
  \left[\begin{array}{ccccc}
   y_1 & y_2  & \cdots & y_{n_2}  \\
   y_2 & y_3  & \cdots & y_{n_2-1} \\
  \vdots & \vdots  & & \vdots \\
   y_{n_1} & y_{n_1-1} & \cdots & y_{n_1+n_2-1}
  \end{array}\right] 
 = \sum\limits_{k=1}^r
 c_k \left[\begin{array}{c}
 1 \\ e^{\j\omega_k} \\ \vdots \\ e^{\j (n_1-1)\omega_k} \end{array}\right]
 \left[\begin{array}{c}
 1 \\ e^{-\j\omega_k} \\ \vdots \\ e^{-\j (n_2-1)\omega_k} 
  \end{array}\right]^H \\*[5ex]
  & |\omega_k| \leq \omega_\mathrm c,\quad k = 1,\ldots,r,
\end{array}
\EEQ
where $n_1+n_2-1=n$.
This problem is equivalent to  
\BEQ \label{e-huber-matrix-sdp}
\begin{array}{ll}
\mbox{minimize} & \gamma\sum\limits_{i=1}^m \phi(y_i - y_{\mathrm m,i})
  + (\Tr V)/(2n_1) + (\Tr W)/(2n_2) \\
\mbox{subject to} & X = \left[\begin{array}{cc}
    V & Y \\ Y^H & W \end{array}\right] \succeq 0 \\*[2ex]
 & F X F^T = GXG^T \\
 & - FXG^T - GXF^T + 2\cos\omega_\mathrm c GXG^T \preceq 0
\end{array}
\EEQ
where $G$ and $F$ are block diagonal with blocks
\[
G_1 = \left[\begin{array}{cc}
I_{n_1-1} & 0
\end{array}\right], \qquad
F_1 = \left[\begin{array}{cc}
0 & I_{n_1-1}
\end{array}\right], \qquad
G_2 = \left[\begin{array}{cc}
0 & I_{n_2-1}
\end{array}\right], \qquad
F_2 = \left[\begin{array}{cc} I_{n_2-1} & 0 \end{array}\right].
\]
The variables in~(\ref{e-huber-matrix-sdp}) are the matrices
$V$, $Y$, $W$.  The elements $y_i$ in the objective are the elements in
the first row and last column of the matrix variable $Y$.
The two SDP~(\ref{e-huber-sdp}) 
and~(\ref{e-huber-matrix-sdp}) give the same result $y$,
but may have different numerical properties (in terms of accuracy or 
complexity).

\subsection{Direction of arrival estimation} \label{s-doa}
This example illustrates the use of frequency interval 
constraints in direction of arrival estimation.
We consider the example described in~\cite[section 3.1]{ChV:16}: 
\BEQ \label{e-doa-bp}
   \begin{array}{ll}
\mbox{minimize} & 
 \sum\limits_{j=1}^3 \sum\limits_{k=1}^{r_j} |x_{jk}| \\*[-.1ex]
\mbox{subject to} 
  &  y_j = \sum\limits_{k=1}^{r_j} x_{jk}\!\left[\begin{array}{c}
   1 \\ e^{\j\pi\sin \theta_{jk}} \\ \vdots \\  e^{\j (n-1)
  \pi \sin \theta_{jk}} \end{array}\right] 
  \\*[4.7ex]
  &  \theta_{jk} \in \Theta_j, \quad k=1,\ldots, r_j, \quad j=1,2,3  
     \\*[.5ex]
   & (y_1 + y_2)_{I_1} = b_1, \quad  (y_2 + y_3)_{I_2} = b_2.
   \end{array}
\EEQ
The vectors $b_1$ and $b_2$ contain the outputs of two subsets
of the elements in a linear array of $n$ non-isotropic antennas. 
Elements in the first group, indexed by the index set $I_1$,
measure input signals arriving from angles in 
$\Theta_1\cup\Theta_2 = [-\pi/2, -\pi/6] \cup [-\pi/6, \pi/6]$.
Elements in the second group, indexed by the index set $I_1$,
measure input signals arriving from 
$\Theta_2\cup\Theta_3=[-\pi/6, \pi/6] \cup [\pi/6, \pi/2]$.
The convex formulation of this problem can be found in
\cite{ChV:16}.

\begin{figure}
\includegraphics[width=.49\linewidth]{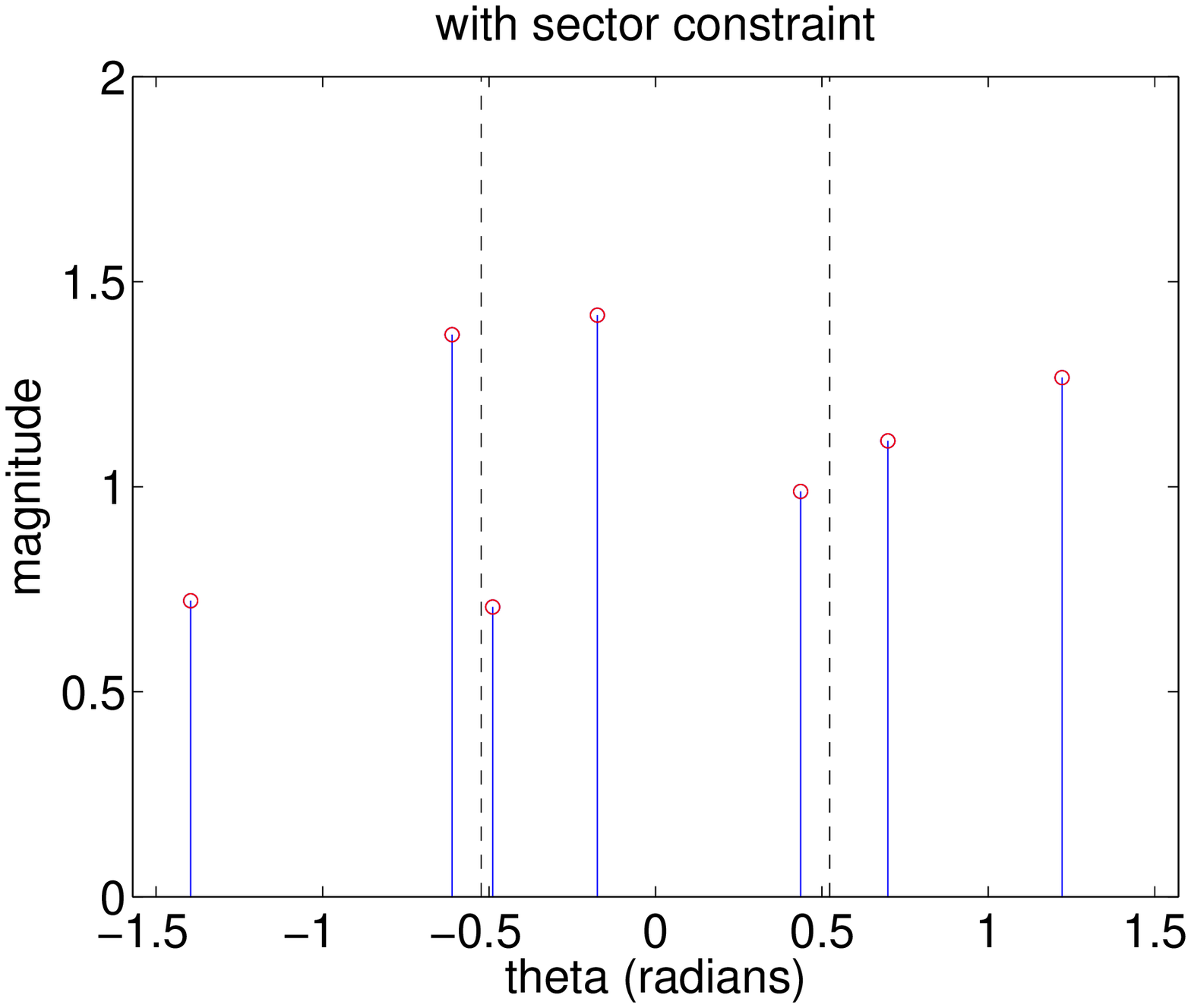}
\hspace*{\fill}
\includegraphics[width=.49\linewidth]{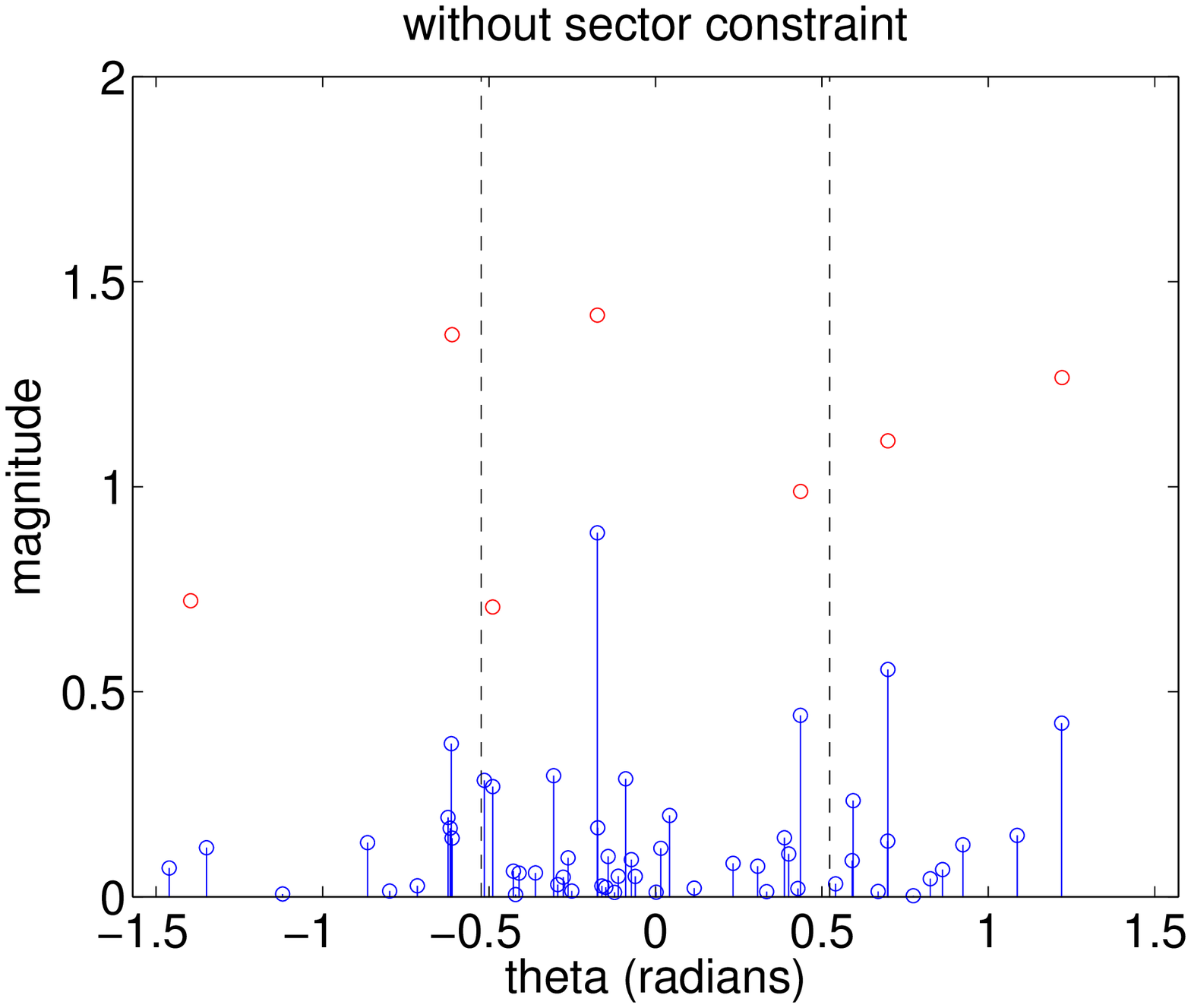}
\caption{Directional of arrival estimation with and without interval
constraints (section~\ref{s-doa}).} \label{f-doa}
\end{figure}

Figure~\ref{f-doa} shows the results of an instance with $n=500$ 
elements in the array,
but using only a total of $40$ randomly selected measurements
($|I_1| = |I_2| = 20$).
The red dots show the angles and magnitudes of 7 signals used
to compute the measurement vectors $b_1$, $b_2$.
The estimated angles and coefficients $|c_{jk}|$ are shown with blue lines.
The right-hand plot shows the solution if we omit the
interval constraints in~(\ref{e-doa-bp}).

\begin{figure}
\centering
\includegraphics[width=.5\linewidth]{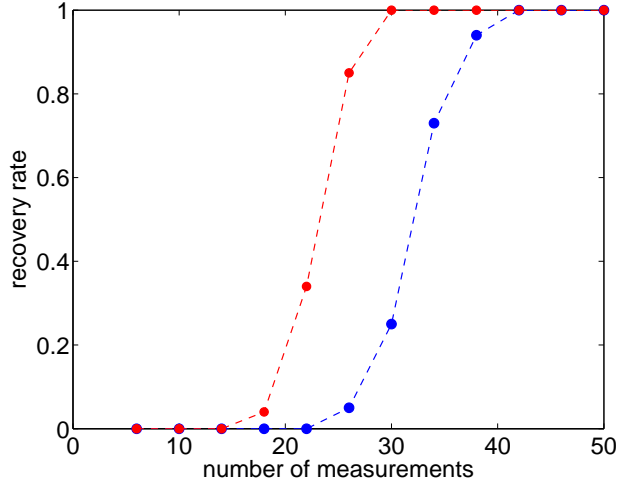}
\caption{Comparison of recovery rate for different number of 
available measurements with interval constraints (red) and without (blue),
in the example of section~\ref{s-doa}.}
\label{f-recovery}
\end{figure}

Figure~\ref{f-recovery} shows the success rate as a function of
the number $|I_1| + |I_2|$ of available measurements, for an example with 
$n=50$ elements, and the same angles as in~\cite{ChV:16} and
figure~\ref{f-doa}. 
Each data point is the average of $100$ trials,
with different, randomly generated coefficients,
and different random selections of the two sensor groups. 
We observe that solving the optimization problem with the interval 
constraints has a higher rate of exact recovery. 
For example, with $30$ available measurements, including the interval 
constraints gave the exact answer in all instances, whereas the method 
without the interval constraints was successful in only about 25\% of 
the instances.

\subsection{Direction of arrival from multiple measurement vectors}
\label{s-doa-multiple}
This example demonstrates the advantage of using multiple measurement 
vectors (or snapshots), as pointed out in~\cite{LC14,YX14}.
Suppose we have $K$ omnidirectional sensors placed at randomly chosen 
positions of a linear grid of length $n$. 
The measurements of the $K$ sensors at one time instance form one 
measurement vector.
We collect $m$ of these measurement vectors, at $m$ different times,
and assume that the directions of arrival and the source magnitudes 
remain constant while the measurements are taken.
The problem is formulated as
\BEQ \label{e-doa-snapshot}
\begin{array}{ll}
\mbox{minimize} & \sum\limits_{k=1}^r \|c_k\| \\
\mbox{subject to} 
 & Y = \sum\limits_{k=1}^r 
 \left[\begin{array}{c}
  1 \\ e^{\j\alpha\sin\theta_k} \\ \vdots \\ 
  e^{\j (n-1)\alpha\sin\theta_k} \end{array} 
 \right] c_k^H \\
 & Y_I = B \\
& |\theta_k| \leq \theta_\mathrm c, \quad k = 1,\ldots,r,
\end{array}
\EEQ
with variables $Y\in\complex^{n\times m}$, $c_k\in\complex^m$, 
$\omega_k$, and $r$.
Here $\alpha = 2\pi d / \lambda_\mathrm c$, where $d$ is the distance between the grid points
and $\lambda_\mathrm c$ is the signal wavelength,
and $\theta_\mathrm c$ is a given cutoff angle.
The columns of the $K\times m$ vector $B$ are the measurement
vectors.  The matrix $Y_I$ is the submatrix of $Y$
containing the rows indexed by $I$.  
The problem can be interpreted as identifying a continuous form
of group sparsity~\cite{Gra15}.
The convex formulation is
\[
\begin{array}{ll}
\mbox{minimize} & (\Tr V)/(2n) + (\Tr W)/2 \\*[1ex]
\mbox{subject to} 
 & \left[\begin{array}{cc}
    V & Y \\ Y^H & W \end{array}\right] \succeq 0 \\*[2ex]
 & FVF^H - GVG^H = 0 \\
 & -FVG^H - GVF^H + 2 \cos \omega_\mathrm c GVG^H \preceq 0 \\
 & Y_I = B 
\end{array}
\]
with $F$ and $G$ defined in~(\ref{e-FG-toep}) and $\omega_\mathrm c = \alpha\sin\theta_\mathrm c$.

\begin{figure}
\hspace*{\fill}
\includegraphics[width=.47\linewidth]{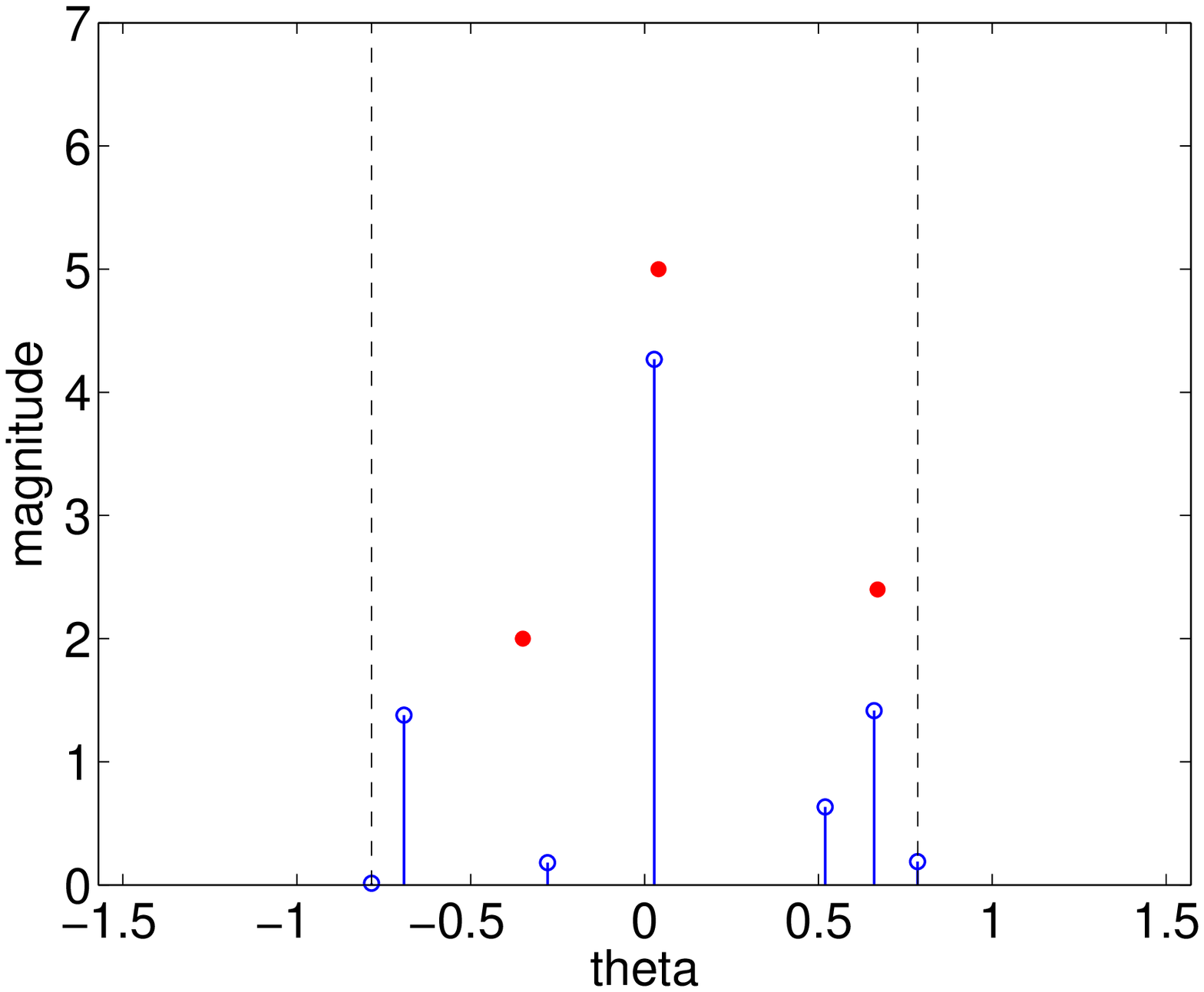}
\hfill
\includegraphics[width=.45\linewidth]{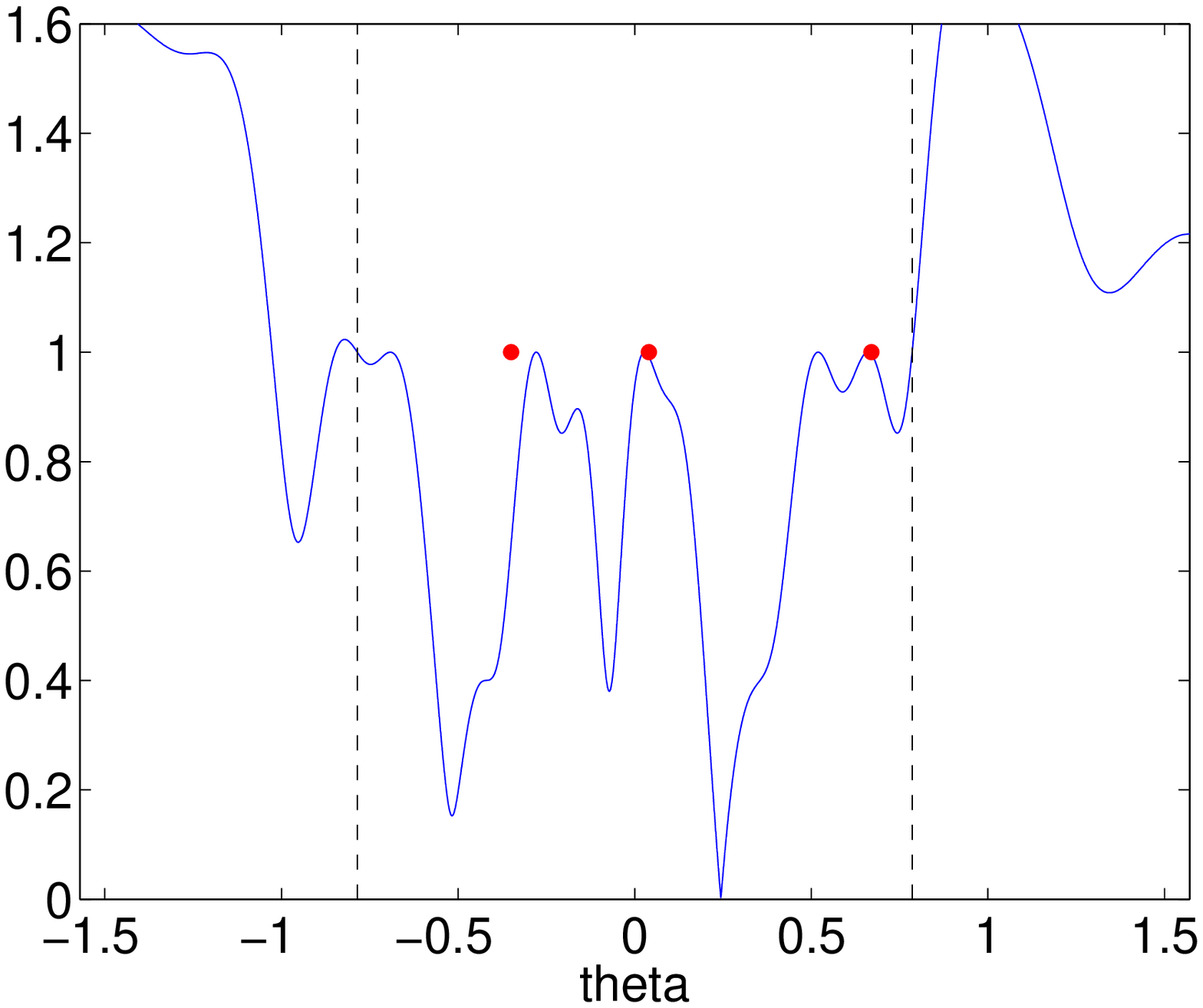} 
\hspace*{\fill}
\\
\hspace*{\fill}
\includegraphics[width=.47\linewidth]{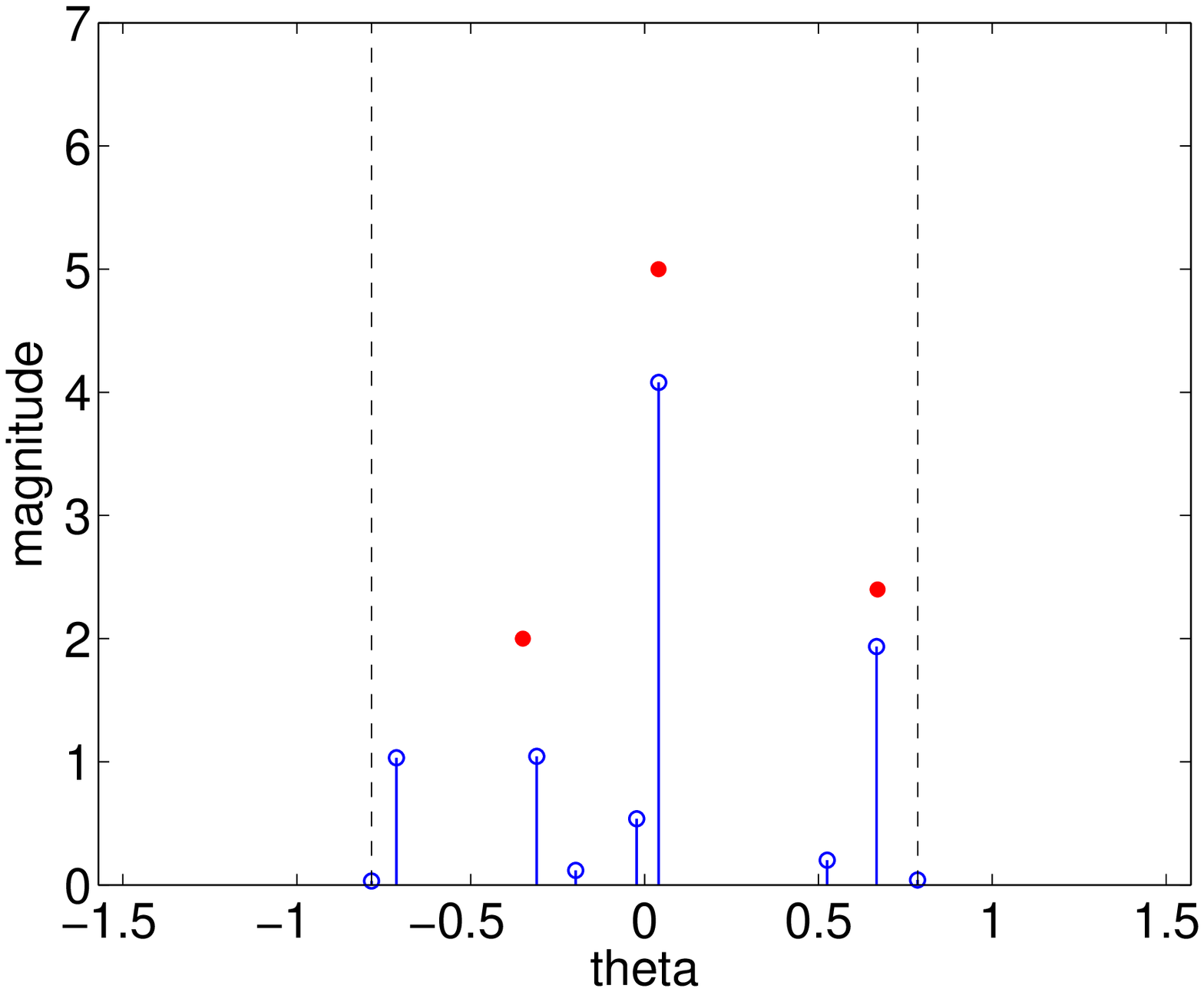}
\hfill
\includegraphics[width=.45\linewidth]{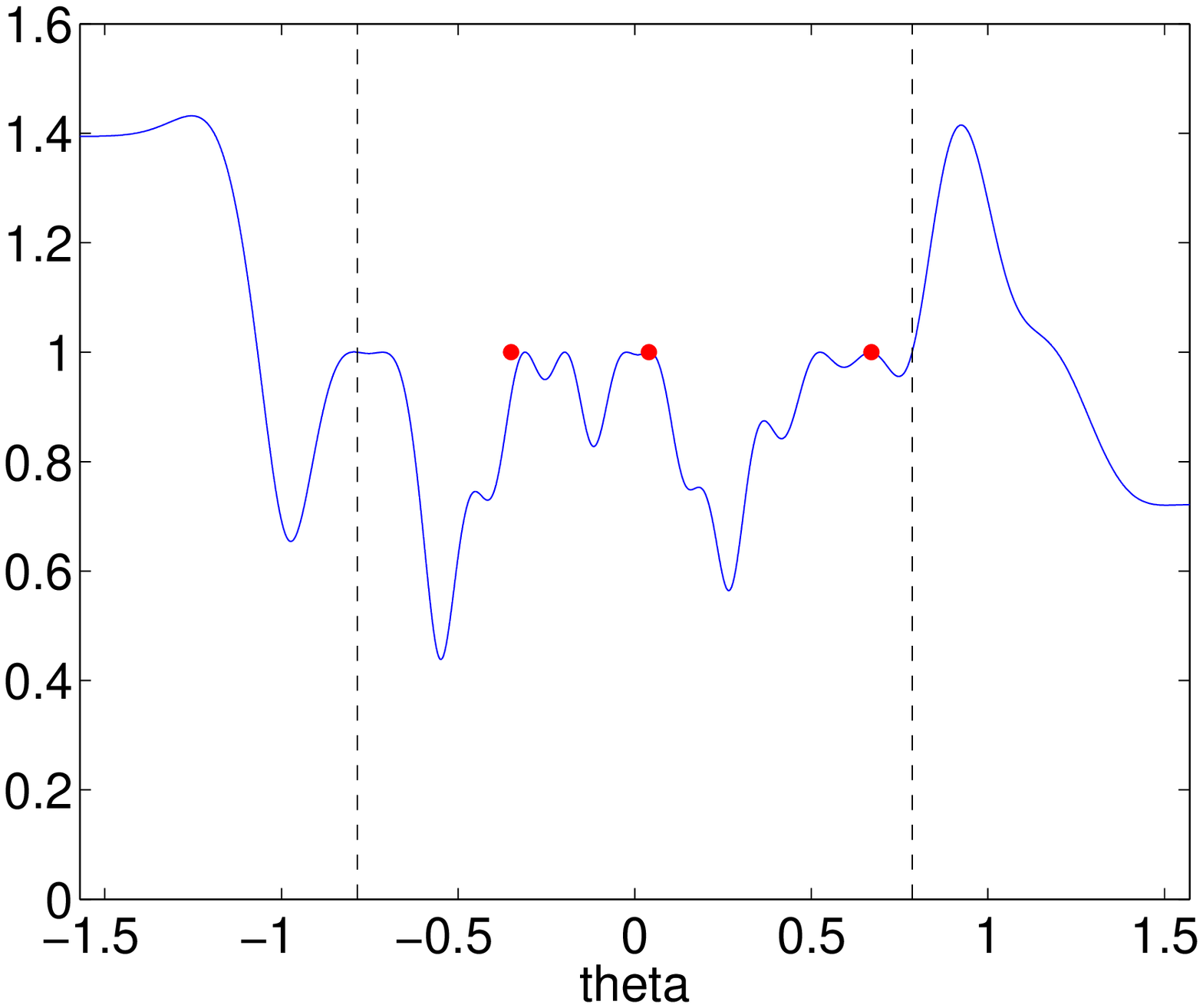} 
\hspace*{\fill}
\\
\hspace*{\fill}
\includegraphics[width=.47\linewidth]{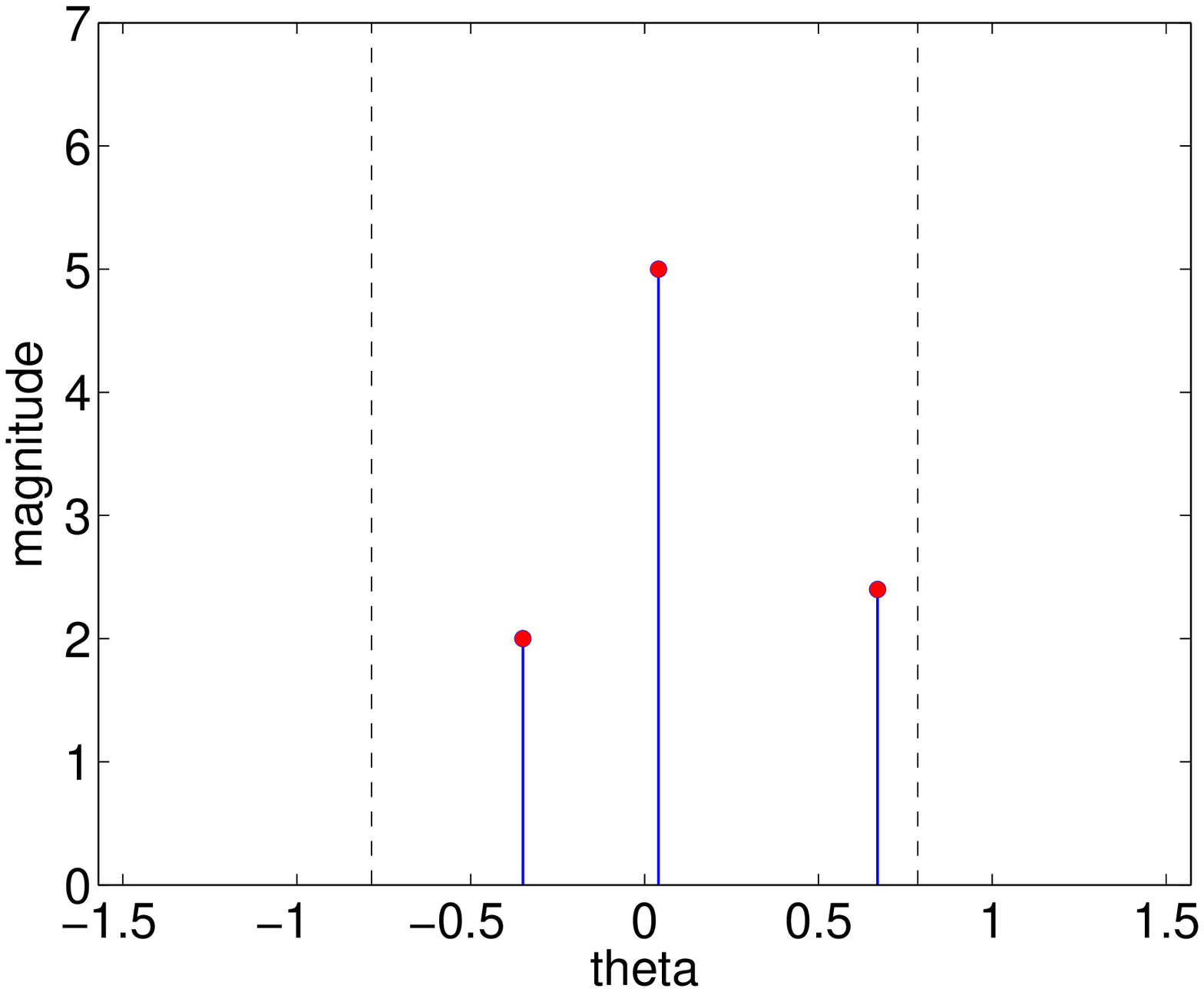}
\hfill
\includegraphics[width=.45\linewidth]{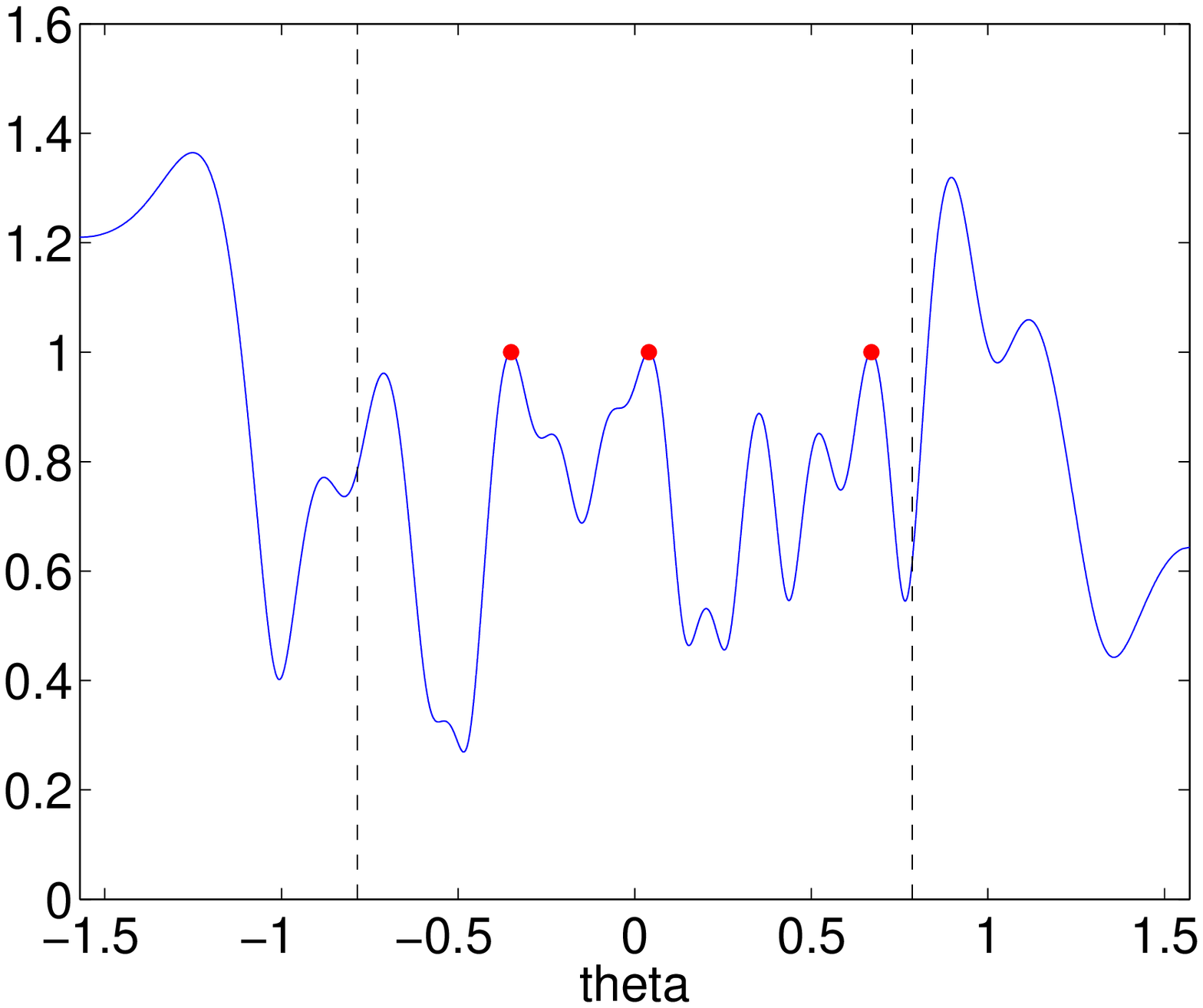}
\hspace*{\fill}
\caption{From top to bottom are shown the results of recovery 
with $1$, $15$, $30$ measurement vectors in the DOA estimation
problem of section~\ref{s-doa-multiple}.
The figures on the right show the magnitude of the
trigonometric polynomials obtained from the dual optimal solutions.
The red dots represent the true directions of arrival (and magnitudes).}
\label{f-doa-snapshots}
\end{figure}

Figure~\ref{f-doa-snapshots} shows 
an example with $n=30$, $K=7$, $\alpha = 2$, 
and $\theta_\mathrm c = \pi/4$.
We show the solution for $m=1$, $m=15$, $m=30$. 
The blue lines show the values of $\omega_k$ and $\|c_k\|/\sqrt m$
computed by solving problem~(\ref{e-doa-snapshot}).

\section{Conclusion}
In this paper we developed semidefinite representations of a class 
of gauge functions and atomic norms for sets
parameterized by linear matrix pencils. 
The formulations extend the semidefinite representation
of the atomic norm associated with the trigonometric moment curve,
which underlies recent results in continuous or `off-the-grid'
compressed sensing.
The main contribution is a self-contained constructive proof of
the semidefinite representations, using techniques developed  
in the literature on the Kalman-Yakubovich-Popov lemma. 
In addition to opening new possible areas of applications in
system theory and control, 
the connection with the KYP lemma is important for numerical
algorithms.  Specialized techniques for solving SDPs derived from
the KYP lemma, for example, by exploiting 
real symmetries and rank-one structure
\cite{GHNV:03,LoP:04,RoV:06,LiV:07,HaV:14}, should be useful in the 
development of fast solvers for the SDPs discussed in this paper.

\appendix

\section{Subsets of the complex plane} \label{s-regions}
In this appendix we explain the notation used in 
equation~(\ref{e-mC}) to describe subsets of the closed complex plane.
Recall that we use the notation
\[
 q_\Theta(\mu,\nu) = \left[\begin{array}{c}
  \mu \\ \nu \end{array}\right]^H
 \left[\begin{array}{cc}
 \Theta_{11} & \Theta_{12} \\ \Theta_{21} & \Theta_{22}
 \end{array}\right] 
 \left[\begin{array}{c} \mu \\ \nu \end{array}\right]
\]
for the quadratic form defined by a Hermitian $2\times 2$ matrix $\Theta$. 

\paragraph{Lines and circles}
If $\Phi$ is a $2\times 2$ Hermitian matrix with $\det \Phi < 0$,
then the  quadratic equation
\BEQ \label{e-quad-eq}
    q_\Phi (\lambda,1) =  0
\EEQ
defines a straight line  (if $\Phi_{11} = 0$) or a circle (if $\Phi_{11}
\neq 0$) in the complex plane.   
Three important special cases are
\[
\Phi_\mathrm u = 
\left[\begin{array}{cc} 1 & 0 \\ 0 & -1 \end{array} \right], \qquad 
\Phi_\mathrm i = 
\left[\begin{array}{cc} 0 & 1 \\ 1 &  0 \end{array} \right], \qquad 
\Phi_\mathrm r = 
\left[\begin{array}{cc} 0 & \j \\ -\j &  0 \end{array} \right], 
\]
for the unit circle, imaginary axis, and real axis, respectively.
Curves defined by two different matrices $\Phi$, $\tilde \Phi$ can
be mapped to one another by applying a
nonsingular congruence transformation $\tilde \Phi = R\Phi R^H$.

When $\Phi_{11} = 0$, we include the point $\lambda = \infty$  in the 
solution set of~(\ref{e-quad-eq}).  Alternatively, one can define
points in the closed complex plane as directions $(\mu,\nu)\neq 0$. 
If $\nu \neq 0$, the pair $(\mu,\nu)$ represents the complex
number $\lambda = \mu/\nu$.  If $\nu=0$, it represents the point
at infinity.   Using this notation, a circle or line in the closed 
complex plane is  defined as the nonzero solution set of a quadratic
equation
\[
q_\Phi (\mu,\nu) =
\left[\begin{array}{c} \mu \\ \nu \end{array}\right]^H
 \Phi
\left[\begin{array}{c} \mu \\ \nu \end{array}\right] =0, 
\]
with $\det\Phi < 0$.
A congruence transformation $\tilde \Phi = R \Phi R^H$ corresponds to
a linear transformation between the sets associated with the 
matrices $\Phi$ and $\tilde\Phi$.

\paragraph{Segments of lines and circles}
The second type of set we encounter is defined by a quadratic equality 
and inequality
\BEQ \label{e-quad-eq-ineq} 
q_\Phi(\lambda,1) = 0, \qquad
q_\Psi(\lambda,1) \leq 0. 
\EEQ
We assume that $\det \Phi < 0$.
If the inequality is redundant (\eg, $\Psi = 0$)
the solution set of~(\ref{e-quad-eq-ineq}) is the line or circle
defined by the equality.
Otherwise it is an arc of a circle, a closed interval of a line, or the 
complement of an open interval of a line.
It includes the point at infinity if $\Phi_{11} = 0$ and 
$\Psi_{11} \leq 0$.
Alternatively, one can use homogeneous coordinates and consider
sets of points $(\mu,\nu)$  that satisfy
\BEQ \label{e-quad-eq-ineq-hom}
 q_\Phi (\mu,\nu) = 0, \qquad
 q_\Psi (\mu,\nu) \leq 0, \qquad
 (\mu,\nu) \neq 0.
\EEQ

For easy reference, we list the most common combinations of 
$\Phi$ and $\Psi$ in tables~\ref{t-phid}--\ref{t-phir} 
\cite{IwH:03,IwH:05}.

\begin{table} 
\begin{center} 
\begin{tabular}{ccc}\toprule
$\angle{\lambda}$ & $\Psi$ & Assumptions \\ 
\midrule
$[a-b, a+b]$ &
$\left[\begin{array}{cc} 0 & -e^{\j a} \\ 
 -e^{-\j a} & 2\cos b\end{array}\right]$ &
$0 \leq b \leq  \pi$ \\
$[a, 2\pi-a]$ &
 $\left[\begin{array}{cc}0 & 1 \\ 1 & -2\cos a\end{array}\right]$  
& $0\leq a\leq\pi$ \\ \bottomrule
\end{tabular}
\end{center}
\caption{Common choices of $\Psi$ with $\Phi = \Phi_\mathrm u$
($\lambda$ on the unit circle).}
\label{t-phid}
\end{table}
\begin{table} 
\begin{center} 
\begin{tabular}{ccc}\toprule
$\Im \lambda$ & $\Psi$ & Assumptions \\ 
\midrule 
$[a, b]$ & 
$\left[\begin{array}{cc} 2 & -\j(a+b)\\ \j(a+b)& 2ab
\end{array}\right]$ & $a\leq b$  \\
$[-\infty, -a]\cup [a,\infty]$&
 $\left[\begin{array}{cc} -1 & 0 \\ 0 & a^2\end{array}\right]$  
 & $a\geq 0$ \\
\bottomrule
\end{tabular}
\end{center}
\caption{Common choices of $\Psi$ with $\Phi = \Phi_\mathrm i$
($\lambda$ imaginary). } \label{t-phic}
\end{table}

\begin{table} 
\begin{center} 
\begin{tabular}{ccc}\toprule
$\lambda$ & $\Psi$ & Assumptions \\
\midrule
$[a, b]$ &
$\left[\begin{array}{cc} 2 & -(a+b)\\ -(a+b) & 2ab
 \end{array}\right]$ &  $a \leq b$ \\
$[-\infty, a]\cup [b,\infty]$&
 $\left[\begin{array}{cc} -2 & a+b \\ a+b & -2ab\end{array}\right]$  
& $a\leq b$
\\ 
$[a, \infty]$ &
 $\left[\begin{array}{cc} 0 & -1 \\ -1 & 2a\end{array}\right]$ \\
$[-\infty, a]$ &
 $\left[\begin{array}{cc} 0 & 1 \\ 1 & -2a\end{array}\right]$  
\\ \bottomrule
\end{tabular}
\end{center}
\caption{Common choices of $\Psi$  with $\Phi = \Phi_\mathrm r$
($\lambda$ real).} \label{t-phir}
\end{table}

As for circles and lines, we can apply a congruence transformation
to reduce~(\ref{e-quad-eq-ineq}) to a simple canonical case.
We mention two examples.  
Iwasaki and Hara \cite[lemma 2]{IwH:05} show 
that for every $\Phi$, $\Psi$ with $\det \Phi < 0$, 
there exists a nonsingular $R$ such that
\BEQ \label{e-IwH-canonical}
 \Phi = R^H \Phi_\mathrm i R, \qquad
 \Psi = R^H \left[\begin{array}{cc} \alpha & \beta \\ \beta
 & \gamma \end{array}\right] R
\EEQ
with $\alpha, \beta, \gamma$ real, and $\alpha \geq \gamma$.
To see this, we first apply a congruence transformation
$\Phi = R_1^H \Phi_\mathrm i R_1$ to 
transform $\Phi$ to $\Phi_\mathrm i$.
Define
\[
 R_1^{-H}\Psi R_1^{-1}
= \left[\begin{array}{cc}
  x & \beta + \j z \\ \beta - \j z & y \end{array}\right]
\]
with real $x$, $y$, $z$, $\beta$, and consider the eigenvalue 
decomposition
\BEQ \label{e-IwH-evd}
 \left[\begin{array}{cc} x & \j z \\ -\j z & y \end{array}\right]
= Q 
\left[ \begin{array}{cc} \alpha & 0 \\ 0 & \gamma \end{array}
\right]Q^H,
\EEQ
with eigenvalues sorted as $\alpha \geq \gamma$.
Since the $2,1$ element of the matrix on the left-hand side
of~(\ref{e-IwH-evd}) is purely imaginary,
the columns of $Q$ can be normalized to be of the form
\[
 Q = \left[\begin{array}{cc}
  u  & \j v \\ 
  \j v  & u \end{array}\right]
\]
with $u$ and $v$ real, and $u^2 + v^2 = 1$.
This implies that $Q \Phi_\mathrm i Q^H = Q^H \Phi_\mathrm i Q =
\Phi_\mathrm i$ and
\[
 Q^H \left[\begin{array}{cc} x & \beta +\j z \\ 
 \beta - \j z & y \end{array}\right]Q
= Q^H \left[\begin{array}{cc} x & \j z \\ 
 - \j z & y \end{array}\right]Q
 + \left[\begin{array}{cc} 
   0 &  \beta \\ \beta & 0 \end{array}\right]
= \left[\begin{array}{cc} \alpha & \beta \\ 
 \beta & \gamma \end{array}\right].
\]
The transformation~(\ref{e-IwH-canonical}) now follows by
taking $R = Q^H R_1$.

Applying the congruence defined by $R$, we can reduce the 
conditions~(\ref{e-quad-eq-ineq-hom}) to an equivalent system
\BEQ \label{e-IwH-transformed-hom}
 \left[\begin{array}{c} \mu' \\ \nu' \end{array}\right]^H
 \left[\begin{array}{cc}  0 & 1 \\ 1 & 0 \end{array}\right]
 \left[\begin{array}{c} \mu' \\ \nu' \end{array}\right]
 = 0, \qquad
 \left[\begin{array}{c} \mu' \\ \nu' \end{array}\right]^H
 \left[\begin{array}{cc}  \alpha & 0 \\ 0 & \gamma \end{array}\right]
 \left[\begin{array}{c} \mu' \\ \nu' \end{array}\right]
 \leq 0, \qquad 
 (\mu', \nu') \neq 0,
\EEQ
where $(\mu', \nu')  = R (\mu,\nu)$. In non-homogeneous coordinates,
\BEQ \label{e-IwH-transformed}
 \Re{\lambda'} = 0, \qquad \alpha |\lambda'|^2 +\gamma \leq 0.
\EEQ
Keeping in mind that $\alpha \geq \gamma$, we can distinguish four
cases.  
If $0 < \gamma \leq \alpha$ the solution set 
of~(\ref{e-IwH-transformed}) is empty.
If $\gamma = 0 < \alpha$ the solution set is a singleton $\{0\}$.
If $\gamma < 0 < \alpha$, the solution set 
of~(\ref{e-IwH-transformed}) is the interval of the imaginary
axis defined by $|\lambda'| \leq (-\gamma/\alpha)^{1/2}$.
If $\gamma \leq \alpha \leq 0$, the inequality is redundant and 
the solution set is the entire imaginary axis.

Another useful canonical form of~(\ref{e-quad-eq-ineq})
is obtained by transforming the solution set to a subset of the unit 
circle.  If we define
\[
 T = \frac{1}{\sqrt 2} \left[\begin{array}{cc}
  1 & 1 \\ -1 & 1 \end{array}\right] R, 
\qquad 
 \epsilon = \frac{1}{2} (\alpha + \gamma),
\qquad 
 \delta = \frac{1}{2} (\alpha - \gamma),
\qquad 
 \eta = \beta.
\]
then it follows from from~(\ref{e-IwH-canonical}) that
\[
 \Phi = T^H \Phi_\mathrm u T, \qquad
 \Psi = T^H \left[\begin{array}{cc} \epsilon +\eta & 
 -\delta \\ -\delta & \epsilon - \eta \end{array}\right] T.
\]
The coefficients $\epsilon$, $\delta$, $\eta$ are real,
with $\delta \geq 0$.
The congruence defined by $T$ therefore transforms
the conditions~(\ref{e-quad-eq-ineq-hom}) to an equivalent system
\[
\left[\begin{array}{c} \mu' \\ \nu' \end{array}\right]^H
\left[\begin{array}{cc} 1 & 0 \\ 0 & -1 \end{array}\right]
\left[\begin{array}{c} \mu' \\ \nu' \end{array}\right] = 0, \qquad
\left[\begin{array}{c} \mu' \\ \nu' \end{array}\right]^H
\left[\begin{array}{cc} 0 & -\delta \\ -\delta & 2\epsilon
 \end{array}\right]
\left[\begin{array}{c} \mu' \\ \nu' \end{array}\right] \leq 0, \qquad
\]
where $(\mu', \nu')  = T (\mu,\nu)$. In non-homogeneous coordinates,
this is
\[
 |\lambda'|^2 = 1, \qquad \delta \Re{\lambda'} \geq 
 \epsilon.
\]
The solution set is empty if $\epsilon > \delta$.
It is the unit circle if $\epsilon \leq -\delta$.
It is the singleton $\{1\}$ if $\epsilon = \delta > 0$.
It is a segment of the unit circle if $-\delta < \epsilon < \delta$.

\section{Matrix factorization results} 
\label{s-matrix-fact}
This appendix contains a self-contained proof of 
Lemma~\ref{l-quad-eq-ineq-general}, needed in the proof
of Theorem~\ref{t-decomp}, and some other matrix factorization 
results that have appeared in papers on the Kalman-Yakubovich-Popov (KYP) 
lemma \cite{Ran:96,IMF:00,BaV:02,BaV:03,PiV:11}.
We include the proofs because their constructive character is
important for the result in Theorem~\ref{t-decomp}.

Lemma~\ref{l-rantzer} is based on 
\cite[lemma 3]{Ran:96} and~\cite[lemma 5]{IwH:05}.
Lemma~\ref{l-quad-eq-ineq-general} can be found in 
\cite[corollary 1]{PiV:11}.
\begin{lemma} \label{l-rantzer}
Let $U$ and $V$ be two matrices in $\complex^{p\times r}$.
\BIT
\item If $UU^H = VV^H$, then $U=V\Lambda$ for some unitary matrix 
 $\Lambda\in\complex^{r\times r}$.  
\item If $UU^H = VV^H$ and $UV^H + VU^H = 0$, 
then $U=V\Lambda$ for some unitary and skew-Hermitian matrix 
$\Lambda\in\complex^{r\times r}$.
\item If $UU^H \preceq VV^H$ and $UV^H + VU^H = 0$, 
then $U=V\Lambda$ for some skew-Hermitian matrix 
$\Lambda\in\complex^{r\times r}$ with $\|\Lambda\|_2 \leq 1$.
\EIT
\end{lemma}

\proof\
If $UU^H = VV^H$, then  $U$ and $V$ have singular value decompositions
of the form
\[ 
  U = P \Sigma Q_u^H, \qquad V = P\Sigma Q_v^H,
\]
with unitary matrices $P\in\complex^{p\times p}$,
diagonal $\Sigma \in \complex^{p\times r}$, 
and unitary $Q_u, Q_v \in \complex^{r\times r}$. 
The unitary matrix $\Lambda = Q_vQ_u^H$ satisfies $U = V\Lambda$.

To show the second part of the lemma, we substitute the singular value 
decompositions of $U$ and $V$ in the equation  $UV^H+V^HU=0$:
\[
 \Sigma (Q_u^H Q_v  + Q_v^H Q_u) \Sigma^T = 0.
\]
We define $\tilde\Lambda = Q_u^HQ_v$ (a unitary $r\times r$ matrix)
and write this  as
\[
\left[\begin{array}{cc} \Sigma_1 & 0 \\ 0 & 0 \end{array}\right]
\left[\begin{array}{cc}
 \tilde \Lambda_{11} + \tilde \Lambda_{11}^H & 
 \tilde \Lambda_{12} + \tilde \Lambda_{21}^H \\ 
 \tilde \Lambda_{21} + \tilde \Lambda_{12}^H & 
 \tilde \Lambda_{22} + \tilde \Lambda_{22}^H 
\end{array}\right]
\left[\begin{array}{cc} \Sigma_1 & 0 \\ 0 & 0 \end{array}\right]
 = 0
\]
with $\Sigma_1$ positive diagonal of size $q\times q$, where 
$q =\Rank(U) = \Rank(V)$,
and $\tilde \Lambda_{11}$ the $q\times q$ leading diagonal block of 
$\tilde \Lambda$.
This shows that $\tilde\Lambda_{11} + \tilde\Lambda_{11}^H =0$, 
so $\tilde \Lambda$ is unitary with a skew-Hermitian $1,1$ block.
Since $\tilde \Lambda_{11}$ is skew-Hermitian it has a Schur decomposition 
$\tilde \Lambda_{11} = Q \Delta Q^H$ 
with unitary $Q\in\complex^{q\times q}$, and 
$\Delta$ a diagonal and purely imaginary matrix.  
Moreover $\Delta\Delta^H \preceq I$ because $\tilde \Lambda_{11}$ is a
submatrix of the unitary matrix $\tilde \Lambda$.
Partition $Q$ and $\Delta$ as
\[
 \tilde \Lambda_{11} = 
 \left[\begin{array}{cc} Q_1 & Q_2 \end{array}\right]
 \left[\begin{array}{cc} \Delta_1 & 0\\ 0 & \Delta_2 \end{array}\right]
 \left[\begin{array}{cc} Q_1 & Q_2 \end{array}\right]^H
\]
with $\Delta_1\Delta_1^H \prec I$ and $\Delta_2 \Delta_2^H = I$.
Since $\tilde \Lambda$ is unitary, we have
\BEAS
 \tilde \Lambda_{12}\tilde \Lambda_{12}^H 
 & = & I - \tilde \Lambda_{11}\tilde \Lambda_{11}^H  \\
 & = & Q_1Q_1^H + Q_2 Q_2^H - Q_1 \Delta_1 \Delta_1^H Q_1^H
    - Q_2 \Delta_2 \Delta_2^H Q_2^H \\
 & = & Q_1(I-\Delta_1\Delta_1^H) Q_1^H,
\EEAS
and, by the first part of the lemma,
$\tilde \Lambda_{12} = Q_1 (I-\Delta_1\Delta_1^H)^{1/2} W$ for
some unitary matrix $W$.
Therefore the matrix
\BEAS
\lefteqn{
\left[\begin{array}{cc}
\tilde \Lambda_{11} & \tilde\Lambda_{12} \\
-\tilde\Lambda_{12}^H & W^H \Delta_1^H W
\end{array}\right] }\\
& = & \left[\begin{array}{ccc} 
  Q_1 & Q_2 & 0 \\ 0 & 0 & W^H \end{array}\right]
  \left[\begin{array}{ccc} 
  \Delta_1 & 0 & (I-\Delta_1\Delta_1^H)^{1/2} \\ 
  0 & \Delta_2 & 0 \\ -(I-\Delta_1\Delta_1^H)^{1/2} & 0 & \Delta_1^H 
 \end{array}\right]
 \left[\begin{array}{cc} Q_1^H & 0 \\ Q_2^H & 0 \\ 0 & W 
 \end{array}\right]
\EEAS
is skew-Hermitian (from the expression on the left-hand side and
the fact that $\tilde \Lambda_{11}$ is skew-Hermitian
and $\Delta_1$ is purely imaginary)
and unitary (the right-hand side is a  product of three unitary matrices).
If we now define
\[
\Lambda = Q_v \left[\begin{array}{cc}
 \tilde \Lambda_{11} & \tilde \Lambda_{12} \\
 -\tilde \Lambda_{12}^H & W^H \Delta_1^H W \end{array}\right] Q_v^H 
\]
then $\Lambda$ is unitary and skew-Hermitian, and 
\BEAS
U 
 & = & P \left[\begin{array}{cc} \Sigma_1 & 0 \\ 0 & 0 \end{array}\right]
  \left[\begin{array}{cc} 
    \tilde \Lambda_{11} & \tilde \Lambda_{12} \\
    \tilde \Lambda_{21} & \tilde \Lambda_{22} \end{array}\right] Q_v^H \\
 & = & P \left[\begin{array}{cc} \Sigma_1 & 0 \\ 0 & 0 \end{array}\right]
  \left[\begin{array}{cc} 
    \tilde \Lambda_{11} & \tilde \Lambda_{12} \\
    -\tilde \Lambda_{12}^H & 
   W^H \Delta_1^H W \end{array}\right] Q_v^H \\
 & = & P \left[\begin{array}{cc} \Sigma_1 & 0 \\ 0 & 0 \end{array}\right]
  Q_v^H \Lambda \\
 & = & V \Lambda.
\EEAS
This proves part two of the lemma.

Assume $UU^H \preceq VV^H$ and $VV^H-UU^H$ has rank $s$.
We use any factorization $VV^H - UU^H = \tilde U\tilde U^H$ with 
$\tilde U\in\complex^{p\times s}$ and 
write $UU^H \preceq VV^H$ and $UV^H+VU^H= 0$ as
\[
\left[\begin{array}{cc} U & \tilde U \end{array}\right]
\left[\begin{array}{cc} U & \tilde U \end{array}\right]^H
= \left[\begin{array}{cc} V & 0 \end{array}\right]
 \left[\begin{array}{cc} V & 0 \end{array}\right]^H
\]
and
\[
\left[\begin{array}{cc} U & \tilde U \end{array}\right]
\left[\begin{array}{cc} V & 0 \end{array}\right]^H
+ 
\left[\begin{array}{cc} V & 0 \end{array}\right]
\left[\begin{array}{cc} U & \tilde U \end{array}\right]^H =  0.
\]
It follows from part 2 that 
\[
\left[\begin{array}{cc} U & \tilde U \end{array}\right]
= \left[\begin{array}{cc} V & 0 \end{array}\right]
 \left[\begin{array}{cc}
  \tilde \Lambda_{11} & \tilde \Lambda_{12} \\ 
  \tilde \Lambda_{21} & \tilde \Lambda_{22} 
 \end{array}\right]
\]
with $\tilde \Lambda$ unitary and skew-Hermitian. The subblock 
$\Lambda= \tilde \Lambda_{11}$ satisfies $U=V\Lambda$, $\Lambda +
\Lambda^H = 0$ and
$\Lambda^H \Lambda \preceq I$.
\qed\medskip

\begin{lemma} \label{l-quad-eq-ineq-general}
Let $\Phi$, $\Psi\in\herm^2$ with $\det\Phi < 0$. 
If $U, V\in\complex^{p\times r}$ satisfy
\BEA
\Phi_{11} UU^H  + \Phi_{21} UV^H + \Phi_{12} VU^H + \Phi_{22} VV^H 
 & = & 0, \label{e-quadeq-1}\\ 
\Psi_{11} UU^H  + \Psi_{21} UV^H + \Psi_{12} VU^H + \Psi_{22} VV^H 
 & \preceq & 0, \label{e-quadeq-2}
\EEA
then there exist a matrix $W\in\complex^{p\times r}$, 
a unitary matrix $Q\in\complex^{r\times r}$, and vectors
$\mu,\nu\in\complex^r$ such that
\BEQ \label{e-quadeq-facts}
 U = W\diag(\mu)Q^H, \qquad V = W\diag(\nu)Q^H, \qquad
\EEQ
and
\BEQ \label{e-quadeq-C}
 q_\Phi(\mu_i,\nu_i) =0, \qquad
 q_\Psi(\mu_i,\nu_i) \leq 0, \qquad
 (\mu_i,\nu_i) \neq 0, \qquad i=1,\ldots, r.
\EEQ
\end{lemma}

\proof\
Suppose $U$ and $V$ are $p\times r$ matrices that
satisfy~(\ref{e-quadeq-1}) and~(\ref{e-quadeq-2}).
As explained in appendix~\ref{s-regions}, there exists a nonsingular $R$ 
such that
\[
 \Phi = R^H \left[\begin{array}{cc} 0 & 1 \\ 1 & 0 \end{array}\right]R,
 \qquad
 \Psi = R^H \left[\begin{array}{cc} \alpha & \beta \\ \beta & \gamma 
  \end{array}\right]R 
\]
with $\beta$ real and $\gamma \leq \alpha$.
Define $S = R_{11} U + R_{12} V$ and $T=R_{21} U + R_{22} V$.
From~(\ref{e-quadeq-1}) and~(\ref{e-quadeq-2}),
\[
\left[\begin{array}{cc} S & T \end{array}\right]
\left[\begin{array}{cc} 0 & I \\ I & 0 \end{array}\right]
\left[\begin{array}{c} S^H \\ T^H \end{array}\right]
= \left[\begin{array}{cc} U & V \end{array} \right] 
 \left[\begin{array}{cc}
 \Phi_{11} I & \Phi_{21} I \\ \Phi_{12} I & \Phi_{22} I \end{array}
\right]
\left[\begin{array}{c} U^H \\ V^H \end{array}\right] = 0
\]
and
\[
\left[\begin{array}{cc} S & T \end{array}\right]
\left[\begin{array}{cc} \alpha I & \beta I \\ \beta I & \gamma I 
\end{array}\right]
\left[\begin{array}{c} S^H \\ T^H \end{array}\right]
= \left[\begin{array}{cc}
 U & V \end{array}\right] \left[\begin{array}{cc}
 \Psi_{11} I & \Psi_{21} I \\ \Psi_{12} I & \Psi_{22} I \end{array}
\right]
\left[\begin{array}{c} U^H \\ V^H \end{array}\right] \preceq 0.
\]
Therefore 
\BEQ \label{e-2eqs-transf}
 ST^H + TS^H = 0, \qquad \alpha SS^H + \gamma TT^H \preceq 0.
\EEQ
We show that this implies that
\BEQ\label{e-2eqs-transf-fact-a}
 S = W\diag(s)Q^H, \qquad T = W\diag(t)Q^H,
\EEQ
for some $W\in\complex^{p\times r}$, unitary $Q\in\complex^{r\times r}$,
and  vectors $s,t\in\complex^{r}$ that satisfy
\BEQ\label{e-2eqs-transf-fact-b}
 s_i \bar t_i + \bar s_i t_i = 0, \qquad
\alpha |s_i|^2 + \gamma |t_i|^2 \leq 0, \qquad
(s_i,t_i) \neq 0, \qquad i=1,\ldots,r.
\EEQ
The result is trivial if $S$ and $T$ are zero, since in that
case we can choose
$W$ zero, and arbitrary $Q$, $s$, $t$.
If at least one of the two matrices is nonzero, then the inequality
in~(\ref{e-2eqs-transf}), combined with $\alpha \geq \gamma$,
implies that $\gamma\leq 0$.
Therefore there are three cases to consider.
\BIT
\item
If $\alpha \leq 0$, we write the equality in~(\ref{e-2eqs-transf})
as
\[
 (S+T)(S+T)^H = (S-T)(S-T)^H.
\]
From Lemma~\ref{l-rantzer}, this implies that 
$S+ T = (S-T) \Lambda$ with $\Lambda$  unitary. Let
$\Lambda = Q\diag(\rho)Q^H$ be the Schur decomposition of $\Lambda$,
with $|\rho_i| = 1$ for $i=1,\ldots, r$.
Define
\[
 W =  (S-T)Q, \qquad s = \frac{1}{2} (\rho + \ones), \qquad
 t = \frac{1}{2} (\rho - \ones). 
\]
\item
If $\gamma = 0 < \alpha$, then $S=0$, and we can take $Q=I$
\[
W=T, \qquad s=0, \qquad t=\ones. 
\]
\item
If $\gamma < 0 < \alpha$, then from Lemma~\ref{l-rantzer},
we have $S = (-\gamma/\alpha)^{1/2} T\Lambda$
for some skew-Hermitian $\Lambda$ with $\Lambda^H\Lambda \preceq I$.
This matrix has a Schur decomposition $\Lambda = Q\diag(\rho) Q^H$ 
with $|\rho_i| \leq 1$ for $j=1,\ldots,r$.
Define 
\[
 W = TQ,  \qquad s = (-\gamma/\alpha)^{1/2} \rho, \qquad t =\ones. 
\]
\EIT
The factorizations of $U$ and $V$ now follow from
\[
\left[\begin{array}{c}
 U \\ V\end{array}\right]
= (R^{-1} \otimes I) \left[\begin{array}{c} S \\ T\end{array}\right]
= (R^{-1} \otimes I) \left[\begin{array}{c} W\diag(s) \\ 
W\diag(t) \end{array}\right]Q^H
=  \left[\begin{array}{c} W\diag(\mu) \\ 
W\diag(\nu) \end{array}\right]Q^H
\]
where $\mu$ and $\nu$ are defined as
\[
\left[\begin{array}{c} \mu_i \\ \nu_i \end{array}\right]
 = R^{-1} \left[\begin{array}{c} s_i \\ t_i \end{array}\right],
 \quad i=1,\ldots,r.
\]
These pairs $(\mu_i,\nu_i)$ are nonzero and satisfy
\[
\left[\begin{array}{c} \mu_i \\ \nu_i \end{array}\right]^H
\Phi
\left[\begin{array}{c} \mu_i \\ \nu_i \end{array}\right]
= \left[\begin{array}{c} s_i \\ t_i \end{array}\right]^H
 \left[\begin{array}{cc} 0 & 1 \\ 1 & 0 \end{array}\right] 
 \left[\begin{array}{c} s_i \\ t_i \end{array}\right] 
= \bar s_i t_i + s_i \bar t_i 
= 0 
\]
and
\[
\left[\begin{array}{c} \mu_i \\ \nu_i \end{array}\right]^H 
\Psi
\left[\begin{array}{c} \mu_i \\ \nu_i \end{array}\right]
= \left[\begin{array}{c} s_i \\ t_i \end{array}\right]^H
 \left[\begin{array}{cc} \alpha & \beta \\ \beta & \gamma 
 \end{array}\right]  
 \left[\begin{array}{c} s_i \\ t_i \end{array}\right] 
= \alpha |s_i|^2 + \beta(\bar s_i t_i + s_i\bar t_i)
 + \gamma |t_i|^2  
\leq 0.
\]
\qed

\section{Strict feasibility}
\label{s-slater}
In this appendix we discuss strict feasibility of the 
constraints $X\succeq 0$, (\ref{e-phi}), (\ref{e-psi}) 
in Theorem~\ref{t-decomp}. 

We assume that the set $\mathcal C$ defined in~(\ref{e-mC}) 
is not empty and not a singleton.
This means that if the inequality $q_\Psi(\mu,\nu) \leq 0$ 
in the definition is not redundant, then there exist points in $\mathcal C$
with $q_\Psi(\mu,\nu) < 0$.
We will distinguish these two cases.
\BIT
\item \emph{Line or circle.}
If the inequality $q_\Psi(\mu,\nu) \leq 0$ in the definition is redundant,
we have
\[
 \mathcal C = \{ (\mu,\nu) \in \complex^2 \mid (\mu,\nu) \neq 0,  
 \; q_\Phi(\mu,\nu) = 0 \},
\]
and $\mathcal C$ is a line or circle in homogeneous coordinates.
In this  case we understand by strict feasibility of $X$ that
\BEQ \label{e-strict-feas-case-1}
X \succ 0, \qquad
\Phi_{11} FXF^H + \Phi_{21} FXG^H + \Phi_{12} GXF^H + \Phi_{22} GXG^H = 0.
\EEQ
We also define $\mathcal C^\circ = \mathcal C$.

\item \emph{Segment of line or circle.}
In the second case, $\mathcal C$ is a proper one-dimensional subset of 
the line or circle defined by $q_\Psi(\mu,\nu) = 0$.  
In this case we define strict feasibility  of $X$ as
\BEQ\label{e-strict-feas-case-2}
\mbox{(\ref{e-strict-feas-case-1})}, \qquad
\Psi_{11} FXF^H + \Psi_{21} FXG^H + \Psi_{12} GXF^H + \Psi_{22} GXG^H 
\prec 0.
\EEQ
We also define 
$\mathcal C^\circ = \{ (\mu,\nu) \neq 0 \mid
  q_\Phi(\mu,\nu) = 0, \; q_\Psi(\mu,\nu) < 0\}$.
\EIT
The conditions on $F$ and $G$ that guarantee strict feasibility
will be expressed in terms of the Kronecker structure of the 
matrix pencil $\lambda G - F$ \cite{Gan:05,Van:79}.
For every matrix pencil there exist nonsingular matrices $P$ and $Q$ 
such that 
\BEA 
\lefteqn{P( \lambda G - F) Q }  \nonumber \\
& = & \left[\begin{array}{ccccccccccc}
   L_{\eta_1}(\lambda)^T & 0 & \cdots & 0 & 0 & 0 & 0 & \cdots & 0 \\
   0 & L_{\eta_2}(\lambda)^T & \cdots & 0 & 0 & 0 & 0 & \cdots & 0 \\
   \vdots & \vdots & \ddots & \vdots & \vdots & \vdots & 
       \vdots & & \vdots \\
   0 & 0 & \cdots & L_{\eta_l}(\lambda)^T & 0 & 0 & 0 & \cdots & 0 \\
   0 & 0 & \cdots & 0  & \lambda B - A & 0 & 0 & \cdots & 0 \\
   0 & 0 & \cdots & 0  & 0 & L_{\epsilon_1}(\lambda) & 0 & \cdots & 0 \\
   0 & 0 & \cdots & 0  & 0 & 0 & L_{\epsilon_2} (\lambda) &\cdots & 0 \\
   \vdots & \vdots & & \vdots  & \vdots & \vdots & \vdots & \ddots & 
   \vdots \\
   0 & 0 & \cdots & 0  & 0 & 0  & 0 & \cdots & L_{\epsilon_r}(\lambda) 
\end{array}\right]
\label{e-kronecker}
\EEA
where 
$L_\epsilon(\lambda)$ is the $\epsilon \times (\epsilon +1)$ pencil
\[
 L_\epsilon(\lambda) = \left[\begin{array}{cccccc}
   \lambda & -1      &  0      &  \cdots &  0 & 0  \\
   0       & \lambda & -1      &  \cdots &  0 & 0 \\
   \vdots  & \vdots  & \vdots  &         &  0 & 0 \\
   0       & 0       & 0       &  \cdots & -1 & 0 \\
   0       & 0       & 0       &  \cdots & \lambda & -1 
 \end{array}\right],
\]
and $\lambda B - A$ is a regular pencil, \ie, it is square and 
$\det(\lambda B - A)$ is not identically zero.
The generalized eigenvalues  of $\lambda B - A$ 
are sometimes referred to as the generalized eigenvalues of the pencil
$\lambda G - F$ \cite[page 16]{IOW:99}.
The parameters $\epsilon_1$, \ldots, $\epsilon_r$ are the
\emph{right Kronecker indices} of the pencil and 
the parameters $\eta_1$, \ldots, $\eta_l$ are the
\emph{left Kronecker indices}.
The \emph{normal rank} of the pencil is equal to $p-l$,
where $p$ is the row dimension of $F$ and $G$.

We show that there exists a strictly feasible  $X$ if and
only if the following two conditions hold.
\begin{enumerate}
\item The normal rank of $\lambda G - F$ is $p$. 
 This means that $l =0$  in~(\ref{e-kronecker}).
\item The generalized eigenvalues  of the pencil $\lambda G - F$ 
 (defined as the generalized eigenvalues of $\lambda B -A$)
 are \emph{nondefective}, \ie, their algebraic multiplicity is equal to 
 the geometric multiplicity, and lie in $\mathcal C^\circ$.
 (More accurately, if $\lambda$ is a finite generalized eigenvalue,
 then $(\lambda, 1) \in \mathcal C^\circ$.  If it is an infinite
 generalized eigenvalue, then $(1,0) \in \mathcal C^\circ$.).
\end{enumerate}
A sufficient but more easily verified condition is that 
$\Rank{(\mu G - \nu F)} = p$ for all $(\mu, \nu) \neq 0$,
\ie, $l=0$ and the block $\lambda B-A$ in~(\ref{e-kronecker})
is not present.

\medskip\proof\
Without loss of generality we can assume that the 
pencil is in the Kronecker canonical form ($P=I$, $Q=I$
in~(\ref{e-kronecker})) and that $\Phi = \Phi_\mathrm u$, so the 
equality constraint in~(\ref{e-strict-feas-case-1}) is 
\BEQ \label{e-FXF=GXG}
 FXF^H = GXG^H.
\EEQ
We first show that the two conditions are necessary.
Assume $X$ is strictly feasible.
Partition $X$ as an $(l+1+r)\times (l+1+r)$  block matrix, with block
dimensions equal to the column dimensions of the $l+1+r$ block columns 
in~(\ref{e-kronecker}).
Suppose $l\geq 1$ and consider the $k$th diagonal block
$X_{kk}$ with $1 \leq k\leq l$.
The $k$th diagonal block of the pencil is
\[
 \lambda G_k - F_k = L_{\eta_k}(\lambda)^T = 
 \lambda \left[\begin{array}{c} I_{\eta_k} \\ 0_{1\times \eta_k} 
 \end{array}\right] -
 \left[\begin{array}{c} 0_{1\times \eta_k} \\ I_{\eta_k} 
 \end{array}\right].
\]
The $k$th diagonal block of~(\ref{e-FXF=GXG}) is
$F_k X_{kk} F_k^H = G_kX_{kk} G_k^H$ or
\[
 \left[\begin{array}{c} 0_{1\times\eta_k} \\ I_{\eta_k} 
 \end{array}\right]
  X_{kk}
 \left[\begin{array}{cc} 0_{\eta_k\times 1} & I_{\eta_k} 
 \end{array}\right]
= 
 \left[\begin{array}{c} I_{\eta_k} \\ 0_{1\times\eta_k} 
 \end{array}\right]
  X_{kk}
 \left[\begin{array}{cc} I_{\eta_k} & 
  0_{\eta_k \times 1} \end{array}\right].
\]
This is impossible since $X_{kk} \succ 0$.  
Hence, if~(\ref{e-FXF=GXG}) holds with $X\succ 0$, then $l=0$.

Next suppose $\det(\mu B - \nu A)=0$ for some $(\mu,\nu) \neq 0$.
If $\nu\neq 0$, then $\mu/\nu$ is a finite generalized eigenvalue
of the pencil $\lambda B-A$; if $\nu=0$ then the pencil has a 
generalized eigenvalue at infinity.
Let $y$ be a corresponding left generalized eigenvector, \ie,
$y^H(\mu B - \nu A) = 0$, while $y^H B$ and $y^HA$ are not both zero
(since $y^HB = y^HA = 0$ would imply that the pencil $\lambda B - A$
is singular).
Define $u^H = y^HB$ if $\nu\neq 0$ and $u^H = y^HA$ otherwise.
This is a nonzero vector.
The first diagonal block of~(\ref{e-FXF=GXG}) is 
\BEQ \label{e-AXA}
AX_{11}A^H  = BX_{11}B^H.
\EEQ
From this it follows that $|\mu|^2 u^H X_{11} u = |\nu|^2 u^HX_{11} u$,
and, since $X_{11} \succ 0$, we have 
$q_\Phi(\mu,\nu) = |\mu|^2 - |\nu|^2 = 0$,
\ie, the generalized eigenvalues are on the unit circle.
In addition, if the inequality in~(\ref{e-strict-feas-case-2})
holds, then
\[
\Psi_{11} AX_{11}A^H + \Psi_{21} AX_{11}B^H + \Psi_{12} BX_{11}A^H 
 + \Psi_{22} BX_{11}B^H \prec 0
\]
and from this, $q_\Psi(\mu,\nu) (u^HX_{11}u) < 0$.
This is only possible if $q_\Psi(\mu,\nu) < 0$.
We conclude that if $\det (\mu B - \nu A) = 0$ for nonzero $(\mu,\nu)$, 
then $(\mu, \nu) \in \mathcal C^\circ$.

Next we show that the generalized eigenvalues of the pencil
$\lambda B - A$ are nondefective.  Since $C^\circ$ is the unit circle
or a subset of the unit circle, there are no infinite generalized
eigenvalues.
Assume the pencil is in Weierstrass canonical form, \ie,
\[
\lambda B-A =
 \left[\begin{array}{cccc}
 (\lambda - \rho_1) I - J_{s_1} & 0 & \cdots & 0 \\
 0 & (\lambda - \rho_2) I - J_{s_2} & \cdots & 0 \\
 \vdots & \vdots & \ddots & \vdots \\
 0 & 0 & \cdots & (\lambda - \rho_t) I - J_{s_t}  
 \end{array}\right],
\]
where $\rho_1$, \ldots, $\rho_t$ are the generalized eigenvalues 
(which satisfy $|\rho_i|=1$), and $J_s$ is the $s\times s$ matrix
\[
 J_s = 
 \left[\begin{array}{cccccc}
   0 & 1 & 0 & \cdots & 0 & 0 \\
   0 & 0 & 1 & \cdots & 0 & 0 \\
   0 & 0 & 0 & \cdots & 0 & 0 \\
   \vdots & \vdots & \vdots & & \vdots & \vdots \\
   0 & 0 & 0 & \cdots & 0 & 1\\
   0 & 0 & 0 & \cdots & 0 & 0
   \end{array} \right].
\]
Then~(\ref{e-AXA}) implies that
\[
 (\rho_i - J_{s_i}) X_{11,i} 
 (\rho_i - J_{s_i})^H  = X_{11,i}
\]
where $X_{11, i}$ is the $i$th diagonal block of $X_{11}$,
if we partition $X_{11}$ as a $t\times t$ block matrix
with $i,j$ block of size of $s_i \times s_j$.
Expanding this gives
\[
 |\rho_i|^2 X_{11,i} - \rho_i X_{11,i} J_{s_i}^T 
 - \bar\rho_i J_{s_i} X_{11,i} + J_{s_i}X_{11,i} J_{s_i}^T = X_{11,i}.
\]
Since $|\rho_i|=1$ this simplifies to
\[
 \rho_i X_{11,i} J_{s_i}^T +
  \bar\rho_i J_{s_i} X_{11,i} = J_{s_i}X_{11,i} J_{s_i}^T.
\]
The last row of the second matrix on the left-hand side and the
last row of the matrix on the right-hand side are zero. 
Therefore the last row of the first matrix on the left  is zero.
However the element in column $s_i-1$ is the last diagonal element
of the positive definite matrix $X_{11,i}$.
Hence, we have a contradiction unless $s_i=1$, \ie, the generalized
eigenvalue $\rho_i$ is nondefective.  
We conclude that the two conditions are necessary.

It remains to show that the conditions are sufficient. 
If the two conditions hold, then 
$\lambda G-F$ has the Kronecker canonical form
\[
\lambda G - F
= \left[\begin{array}{cccccc}
 \lambda - \rho_1 & \cdots & 0 & 0 & \cdots & 0 \\  
 \vdots & \ddots & \vdots & \vdots & & \vdots \\
 0 & \cdots & \lambda - \rho_t & 0 & \cdots & 0 \\
 0 & \cdots & 0 & L_{\epsilon_1}(\lambda) & \cdots & 0 \\
 \vdots & & \vdots  & \vdots & \ddots & \vdots  \\ 
 0 & \cdots & 0 & 0 & \cdots & L_{\epsilon_r}(\lambda)  
\end{array}\right]
\]
with $\rho_i\in\mathcal C^\circ$ for $i=1,\ldots, t$.
Define a block diagonal matrix 
\[
X = 
\left[\begin{array}{cccccc}
 1 & \cdots & 0 & 0 & \cdots & 0 \\
 \vdots & \ddots & \vdots & \vdots & & \vdots \\
 0 & \cdots &  1 & 0 & \cdots & 0 \\
 0 & \cdots & 0 & X_{11} & \cdots & 0 \\
 \vdots & & \vdots & \vdots & \ddots & \vdots \\
 0 & \cdots & 0 & 0 & \cdots & X_{rr} 
\end{array}\right]
\]
with
diagonal blocks
\[
 X_{kk} = \sum_{i=1}^{\epsilon_k+1}
  \left[\begin{array}{c} 
  1 \\ \lambda_{ki} \\ \lambda_{ki}^2 \\ \vdots \\
  \lambda_{ki}^{\epsilon_k} \end{array}\right]
  \left[\begin{array}{c} 
  1 \\ \lambda_{ki} \\ \lambda_{ki}^2 \\ \vdots \\
  \lambda_{ki}^{\epsilon_k} \end{array}\right]^H
\]
for $k=1,\ldots, r$, where $\lambda_{k1}$, \ldots,
$\lambda_{k,\epsilon_k+1}$ are distinct elements of $\mathcal C^\circ$.
This matrix $X$ is strictly feasible.
\qed

\end{document}